\documentclass[12pt,twoside]{amsart}
\usepackage{amsmath, amsthm, amscd, amsfonts, amssymb, graphicx}
\usepackage[bookmarksnumbered, plainpages]{hyperref}

\newtheorem{thm}{Theorem}[section]

\newtheorem{lem}[thm]{Lemma}

\newtheorem{defn}[thm]{Definition}

\numberwithin{equation}{section}

\begin{document}

\title{\bf Canonical connections and algebraic Ricci solitons of three-dimensional Lorentzian Lie groups}
\author{Yong Wang}

\thanks{{\scriptsize
\hskip -0.4 true cm \textit{2010 Mathematics Subject Classification:}
53C40; 53C42.
 }}

\maketitle

\begin{abstract}
 In this paper, we compute canonical connections and Kobayashi-Nomizu connections and their curvature on three-dimensional Lorentzian Lie groups with
 some product structure. We define algebraic Ricci solitons associated to canonical connections and Kobayashi-Nomizu connections. We classify algebraic Ricci solitons associated to canonical connections and Kobayashi-Nomizu connections on three-dimensional Lorentzian Lie groups with
 some product structure.
\end{abstract}

\vskip 0.2 true cm


\pagestyle{myheadings}
\markboth{\rightline {\scriptsize Wang}}
         {\leftline{\scriptsize Canonical connections and algebraic Ricci solitons}}

\bigskip
\bigskip


\section{ Introduction}
\indent The concept of the algebraic Ricci soliton was first introduced by Lauret in Riemannian case in \cite{La}. In \cite{La}, Lauret studied the relation
between solvsolitons and Ricci solitons on Riemannian manifolds. More precisely, he proved that any Riemannian solvsoliton metric is a Ricci soliton.
The concept of the algebraic Ricci soliton was extended to the pseudo-Riemannian case in \cite{BO}. In \cite{BO}, Batat and Onda studied
algebraic Ricci solitons of three-dimensional Lorentzian Lie groups. They got a complete classification of algebraic Ricci solitons of three-dimensional Lorentzian Lie groups and they proved that, contrary to the Riemannian case, Lorentzian Ricci solitons needed not be algebraic Ricci solitons.
In \cite{On}, Onda provided a study of algebraic Ricci solitons in the pseudo-Riemannian case and obtained a steady algebraic Ricci soliton and a shrinking algebraic Ricci soliton in the Lorentzian setting. In \cite{ES}, Etayo and Santamaria studied some affine connections on manifolds with the product structure
or the complex structure. In particular, the canonical connection and the Kobayashi-Nomizu connection for a product structure was studied. In this paper,
we introduce a particular product structure on three-dimensional Lorentzian Lie groups and we compute canonical connections and Kobayashi-Nomizu connections and their curvature on three-dimensional Lorentzian Lie groups with
  this product structure. We define algebraic Ricci solitons associated to canonical connections and Kobayashi-Nomizu connections. We classify algebraic Ricci solitons associated to canonical connections and Kobayashi-Nomizu connections on three-dimensional Lorentzian Lie groups with
 this product structure.\\
\indent In Section 2, We classify algebraic Ricci solitons associated to canonical connections and Kobayashi-Nomizu connections on three-dimensional unimodular Lorentzian Lie groups with
 the product structure.
In Section 3, we classify algebraic Ricci solitons associated to canonical connections and Kobayashi-Nomizu connections on three-dimensional non-unimodular Lorentzian Lie groups with
 the product structure.


\vskip 1 true cm

\section{ Algebraic Ricci solitons associated to canonical connections and Kobayashi-Nomizu connections on three-dimensional unimodular Lorentzian Lie groups}

Three-dimensional Lorentzian Lie groups had been classified in \cite{Ca1,CP}(see Theorem 2.1 and Theorem 2.2 in \cite{BO}). Throughout this paper,
we shall by $\{G_i\}_{i=1,\cdots,7}$, denote the connected, simply connected three-dimensional Lie group equipped with a left-invariant Lorentzian metric $g$ and
having Lie algebra $\{\mathfrak{g}\}_{i=1,\cdots,7}$. Let $\nabla$ be the Levi-Civita connection of $G_i$ and $R$ its curvature tensor, taken with the convention
\begin{equation}
R(X,Y)Z=\nabla_X\nabla_YZ-\nabla_Y\nabla_XZ-\nabla_{[X,Y]}Z.
\end{equation}
The Ricci tensor of $(G_i,g)$ is defined by
\begin{equation}\rho(X,Y)=-g(R(X,e_1)Y,e_1)-g(R(X,e_2)Y,e_2)+g(R(X,e_3)Y,e_3),
\end{equation}
where $\{e_1,e_2,e_3\}$ is a pseudo-orthonormal basis, with $e_3$ timelike and the Ricci operator Ric is given by
\begin{equation}\rho(X,Y)=g({\rm Ric}(X),Y).
\end{equation}
We define a product structure $J$ on $G_i$ by
\begin{equation}Je_1=e_1,~Je_2=e_2,~Je_3=-e_3,
\end{equation}
then $J^2={\rm id}$ and $g(Je_j,Je_j)=g(e_j,e_j)$. By \cite{ES}, we define the canonical connection and the Kobayashi-Nomizu connection as follows:
\begin{equation}\nabla^0_XY=\nabla_XY-\frac{1}{2}(\nabla_XJ)JY,
\end{equation}
\begin{equation}\nabla^{{\rm 1}}_XY=\nabla^0_XY-\frac{1}{4}[(\nabla_YJ)JX-(\nabla_{JY}J)X].
\end{equation}
We define\begin{equation}
R^0(X,Y)Z=\nabla^0_X\nabla^0_YZ-\nabla^0_Y\nabla^0_XZ-\nabla^0_{[X,Y]}Z,
\end{equation}
\begin{equation}
R^{{1}}(X,Y)Z=\nabla^{{ 1}}_X\nabla^1_YZ-\nabla^1_Y\nabla^1_XZ-\nabla^1_{[X,Y]}Z.
\end{equation}
The Ricci tensors of $(G_i,g)$ associated to the canonical connection and the Kobayashi-Nomizu connection
are defined by
\begin{equation}\rho^0(X,Y)=-g(R^0(X,e_1)Y,e_1)-g(R^0(X,e_2)Y,e_2)+g(R^0(X,e_3)Y,e_3),
\end{equation}
\begin{equation}\rho^1(X,Y)=-g(R^1(X,e_1)Y,e_1)-g(R^1(X,e_2)Y,e_2)+g(R^1(X,e_3)Y,e_3).
\end{equation}
 The Ricci operators ${\rm Ric}^0$ and ${\rm Ric}^1$ is given by
\begin{equation}\rho^0(X,Y)=g({\rm Ric}^0(X),Y),~~\rho^1(X,Y)=g({\rm Ric}^1(X),Y).
\end{equation}
Let
\begin{equation}\widetilde{\rho}^0(X,Y)=\frac{{\rho}^0(X,Y)+{\rho}^0(Y,X)}{2},~~\widetilde{\rho}^1(X,Y)=\frac{{\rho}^1(X,Y)+{\rho}^1(Y,X)}{2},
\end{equation}
and
\begin{equation}\widetilde{\rho}^0(X,Y)=g(\widetilde{{\rm Ric}}^0(X),Y),~~\widetilde{\rho}^1(X,Y)=g(\widetilde{{\rm Ric}}^1(X),Y).
\end{equation}

\begin{defn}
$(G_i,g,J)$ is called the first (resp. second) kind algebraic Ricci soliton associated to the connection $\nabla^0$ if it satisfies
\begin{equation}
{\rm Ric}^0=c{\rm Id}+D~~({\rm resp.}~\widetilde{{\rm Ric}}^0=c{\rm Id}+D),
\end{equation}
where $c$ is a real number, and $D$ is a derivation of $\mathfrak{g}$, that is
\begin{equation}D[X,Y]=[DX,Y]+[X,DY]~~{\rm for }~ X,Y\in \mathfrak{g}.
\end{equation}
$(G_i,g,J)$ is called the first (resp. second) kind algebraic Ricci soliton associated to the connection $\nabla^1$ if it satisfies
\begin{equation}
{\rm Ric}^1=c{\rm Id}+D~~({\rm resp.}~\widetilde{{\rm Ric}}^1=c{\rm Id}+D).
\end{equation}
\end{defn}
\vskip 0.5 true cm
\noindent{\bf 2.1 Algebraic Ricci solitons of $G_1$}\\
\vskip 0.5 true cm
By (2.1) and Lemma 3.1 in \cite{BO}, we have for $G_1$, there exists a pseudo-orthonormal basis $\{e_1,e_2,e_3\}$ with $e_3$ timelike such that the Lie
algebra of $G_1$ satisfies
\begin{equation}
[e_1,e_2]=\alpha e_1-\beta e_3,~~[e_1,e_3]=-\alpha e_1-\beta e_2,~~[e_2,e_3]=\beta e_1+\alpha e_2+\alpha e_3,~~\alpha\neq 0.
\end{equation}
\vskip 0.5 true cm
\begin{lem}(\cite{Ca},\cite{BO})
The Levi-Civita connection $\nabla$ of $G_1$ is given by
\begin{align}
&\nabla_{e_1}e_1=-\alpha e_2-\alpha e_3,~~\nabla_{e_2}e_1=\frac{\beta}{2}e_3,~~\nabla_{e_3}e_1=\frac{\beta}{2}e_2,\\\notag
&\nabla_{e_1}e_2=\alpha e_1-\frac{\beta}{2} e_3,~~\nabla_{e_2}e_2=\alpha e_3,~~\nabla_{e_3}e_2=-\frac{\beta}{2}e_1-\alpha e_3,\\\notag
&\nabla_{e_1}e_3=-\alpha e_1-\frac{\beta}{2} e_2,~~\nabla_{e_2}e_3=\frac{\beta}{2}e_1+\alpha e_2,~~\nabla_{e_3}e_3=-\alpha e_2.
\notag
\end{align}
\end{lem}
By (2.4) and (2.18), we have
\vskip 0.5 true cm
\begin{lem}
For $G_1$, the following equalities hold
\begin{align}
&\nabla_{e_1}(J)e_1=-2\alpha e_3,~~\nabla_{e_1}(J)e_2=-\beta e_3,~~\nabla_{e_1}(J)e_3=2\alpha e_1+\beta e_2,~~\\\notag
&\nabla_{e_2}(J)e_1=\beta e_3,~~\nabla_{e_2}(J)e_2=2\alpha e_3,~~\nabla_{e_2}(J)e_3=-\beta e_1-2\alpha e_2,~~\\\notag
&\nabla_{e_3}(J)e_1=0,~~\nabla_{e_3}(J)e_2=-2\alpha e_3,~~\nabla_{e_3}(J)e_3=2\alpha e_2.~~\notag
\end{align}
\end{lem}
By (2.4),(2.5), Lemma 2.2 and Lemma 2.3, we have
\vskip 0.5 true cm
\begin{lem}
The canonical connection $\nabla^0$ of $(G_1,J)$ is given by
\begin{align}
&\nabla^0_{e_1}e_1=-\alpha e_2,~~\nabla^0_{e_1}e_2=\alpha e_1,~~\nabla^0_{e_1}e_3=0,\\\notag
&\nabla^0_{e_2}e_1=0,~~\nabla^0_{e_2}e_2=0,~~\nabla^0_{e_2}e_3=0,\\\notag
&\nabla^0_{e_3}e_1=\frac{\beta}{2}e_2,~~\nabla^0_{e_3}e_2=-\frac{\beta}{2}e_1,~~\nabla^0_{e_3}e_3=0.
\notag
\end{align}
\end{lem}
By (2.7) and Lemma 2.4, we have
\vskip 0.5 true cm
\begin{lem}
The curvature $R^0$ of the canonical connection $\nabla^0$ of $(G_1,J)$ is given by
\begin{align}
&R^0(e_1,e_2)e_1=(\alpha^2+\frac{\beta^2}{2})e_2,~~R^0(e_1,e_2)e_2=-(\alpha^2+\frac{\beta^2}{2})e_1,~~R^0(e_1,e_2)e_3=0,\\\notag
&R^0(e_1,e_3)e_1=-\alpha^2 e_2,~~R^0(e_1,e_3)e_2=\alpha^2 e_1,~~R^0(e_1,e_3)e_3=0,\\\notag
&R^0(e_2,e_3)e_1=\frac{\alpha\beta}{2}e_2,~~R^0(e_2,e_3)e_2=-\frac{\alpha\beta}{2}e_1,~~R^0(e_2,e_3)e_3=0.\notag
\notag
\end{align}
\end{lem}
By (2.9),(2.11) and Lemma 2.5, we get
\begin{align}
{\rm Ric}^0\left(\begin{array}{c}
e_1\\
e_2\\
e_3
\end{array}\right)=\left(\begin{array}{ccc}
-\left(\alpha^2+\frac{\beta^2}{2}\right)&0&0\\
0&-\left(\alpha^2+\frac{\beta^2}{2}\right)&0\\
\frac{\alpha\beta}{2}&\alpha^2&0
\end{array}\right)\left(\begin{array}{c}
e_1\\
e_2\\
e_3
\end{array}\right).
\end{align}
If $(G_1,g,J)$ is the first kind algebraic Ricci soliton associated to the connection $\nabla^0$, then
${\rm Ric}^0=c{\rm Id}+D$, so
\begin{align}
\left\{\begin{array}{l}
De_1=-(\alpha^2+\frac{\beta^2}{2}+c)e_1,\\
De_2=-(\alpha^2+\frac{\beta^2}{2}+c)e_2,\\
De_3=\frac{\alpha\beta}{2}e_1+\alpha^2e_2-ce_3.\\
\end{array}\right.
\end{align}
By (2.15) and (2.23), we get $\alpha^2+c=0$, $\beta=0$. Then
\vskip 0.5 true cm
\begin{thm}
$(G_1,g,J)$ is the first kind algebraic Ricci soliton associated to the connection $\nabla^0$ if and only if  $\alpha^2+c=0$, $\beta=0$, $\alpha\neq 0$.
In particular,
\begin{align}
{\rm Ric}^0\left(\begin{array}{c}
e_1\\
e_2\\
e_3
\end{array}\right)=\left(\begin{array}{ccc}
-\alpha^2&0&0\\
0&-\alpha^2&0\\
0&\alpha^2&0
\end{array}\right)\left(\begin{array}{c}
e_1\\
e_2\\
e_3
\end{array}\right),~~~
{D}\left(\begin{array}{c}
e_1\\
e_2\\
e_3
\end{array}\right)=\left(\begin{array}{ccc}
0&0&0\\
0&0&0\\
0&\alpha^2&\alpha^2
\end{array}\right)\left(\begin{array}{c}
e_1\\
e_2\\
e_3
\end{array}\right).
\end{align}
\end{thm}
By (2.12),(2.13) and (2.22), we have
\begin{align}
\widetilde{{\rm Ric}}^0\left(\begin{array}{c}
e_1\\
e_2\\
e_3
\end{array}\right)=\left(\begin{array}{ccc}
-\left(\alpha^2+\frac{\beta^2}{2}\right)&0&-\frac{\alpha\beta}{4}\\
0&-\left(\alpha^2+\frac{\beta^2}{2}\right)&-\frac{\alpha^2}{2}\\
\frac{\alpha\beta}{4}&\frac{\alpha^2}{2}&0
\end{array}\right)\left(\begin{array}{c}
e_1\\
e_2\\
e_3
\end{array}\right).
\end{align}
If $(G_1,g,J)$ is the second kind algebraic Ricci soliton associated to the connection $\nabla^0$, then
$\widetilde{{\rm Ric}}^0=c{\rm Id}+D$, so
\begin{align}
\left\{\begin{array}{l}
De_1=-(\alpha^2+\frac{\beta^2}{2}+c)e_1-\frac{\alpha\beta}{4}e_3,\\
De_2=-(\alpha^2+\frac{\beta^2}{2}+c)e_2-\frac{\alpha^2}{2}e_3,\\
De_3=\frac{\alpha\beta}{4}e_1+\frac{\alpha^2}{2}e_2-ce_3.\\
\end{array}\right.
\end{align}
By (2.15) and (2.26), we get $\frac{\alpha^2}{2}+c=0$, $\beta=0$. Then
\vskip 0.5 true cm
\begin{thm}
$(G_1,g,J)$ is the second kind algebraic Ricci soliton associated to the connection $\nabla^0$ if and only if  $\frac{\alpha^2}{2}+c=0$, $\beta=0$, $\alpha\neq 0$.
In particular,
\begin{align}
\widetilde{{\rm Ric}}^0\left(\begin{array}{c}
e_1\\
e_2\\
e_3
\end{array}\right)=\left(\begin{array}{ccc}
-\alpha^2&0&0\\
0&-\alpha^2&-\frac{\alpha^2}{2}\\
0&\frac{\alpha^2}{2}&0
\end{array}\right)\left(\begin{array}{c}
e_1\\
e_2\\
e_3
\end{array}\right),~~~
{D}\left(\begin{array}{c}
e_1\\
e_2\\
e_3
\end{array}\right)=\left(\begin{array}{ccc}
-\frac{\alpha^2}{2}&0&0\\
0&-\frac{\alpha^2}{2}&-\frac{\alpha^2}{2}\\
0&\frac{\alpha^2}{2}&\frac{\alpha^2}{2}
\end{array}\right)\left(\begin{array}{c}
e_1\\
e_2\\
e_3
\end{array}\right).
\end{align}
\end{thm}

\indent By (2.6), Lemma 2.3 and Lemma 2.4, we have
\vskip 0.5 true cm
\begin{lem}
The Kobayashi-Nomizu connection $\nabla^1$ of $(G_1,J)$ is given by
\begin{align}
&\nabla^1_{e_1}e_1=-\alpha e_2,~~\nabla^1_{e_1}e_2=\alpha e_1,~~\nabla^1_{e_1}e_3=0,\\\notag
&\nabla^1_{e_2}e_1=0,~~\nabla^1_{e_2}e_2=0,~~\nabla^1_{e_2}e_3=\alpha e_3,\\\notag
&\nabla^0_{e_3}e_1=\alpha e_1+\beta e_2,~~\nabla^1_{e_3}e_2=-\alpha e_2-\beta e_1,~~\nabla^1_{e_3}e_3=0.
\notag
\end{align}
\end{lem}
By (2.8) and Lemma 2.8, we have
\vskip 0.5 true cm
\begin{lem}
The curvature $R^1$ of the Kobayashi-Nomizu connection $\nabla^1$ of $(G_1,J)$ is given by
\begin{align}
&R^1(e_1,e_2)e_1=\alpha\beta e_1+(\alpha^2+{\beta^2})e_2,~~R^1(e_1,e_2)e_2=-(\alpha^2+{\beta^2})e_1-\alpha\beta e_2,~~R^1(e_1,e_2)e_3=0,\\\notag
&R^1(e_1,e_3)e_1=-3\alpha^2 e_2,~~R^1(e_1,e_3)e_2=-\alpha^2 e_1,~~R^1(e_1,e_3)e_3=\alpha\beta e_3,\\\notag
&R^1(e_2,e_3)e_1=-\alpha^2 e_1,~~R^1(e_2,e_3)e_2=\alpha^2 e_2,~~R^1(e_2,e_3)e_3=-\alpha^2 e_3.\notag
\notag
\end{align}
\end{lem}
By (2.10),(2.11) and Lemma 2.9, we get
\begin{align}
{\rm Ric}^1\left(\begin{array}{c}
e_1\\
e_2\\
e_3
\end{array}\right)=\left(\begin{array}{ccc}
-\left(\alpha^2+{\beta^2}\right)&\alpha\beta&\alpha\beta\\
\alpha\beta&-\left(\alpha^2+{\beta^2}\right)&-\alpha^2\\
0&0&0
\end{array}\right)\left(\begin{array}{c}
e_1\\
e_2\\
e_3
\end{array}\right).
\end{align}
If $(G_1,g,J)$ is the first kind algebraic Ricci soliton associated to the connection $\nabla^1$, then
${\rm Ric}^1=c{\rm Id}+D$, so

\begin{align}
\left\{\begin{array}{l}
De_1=-(\alpha^2+{\beta^2}+c)e_1+\alpha\beta e_2+\alpha\beta e_3,\\
De_2=\alpha\beta e_1-(\alpha^2+{\beta^2}+c)e_2-\alpha^2 e_3,\\
De_3=-ce_3.\\
\end{array}\right.
\end{align}
By (2.15) and (2.31), we get $\beta=0$, $c=0$. Then
\vskip 0.5 true cm
\begin{thm}
$(G_1,g,J)$ is the first kind algebraic Ricci soliton associated to the connection $\nabla^1$ if and only if $\beta=0$, $c=0$, $\alpha\neq 0$.
In particular,
\begin{align}
{\rm Ric}^1\left(\begin{array}{c}
e_1\\
e_2\\
e_3
\end{array}\right)=\left(\begin{array}{ccc}
-\alpha^2&0&0\\
0&-\alpha^2&-\alpha^2\\
0&0&0
\end{array}\right)\left(\begin{array}{c}
e_1\\
e_2\\
e_3
\end{array}\right),~~~
{D}\left(\begin{array}{c}
e_1\\
e_2\\
e_3
\end{array}\right)=\left(\begin{array}{ccc}
-\alpha^2&0&0\\
0&-\alpha^2&-\alpha^2\\
0&0&0
\end{array}\right)\left(\begin{array}{c}
e_1\\
e_2\\
e_3
\end{array}\right).
\end{align}
\end{thm}
\vskip 0.5 true cm

\indent By (2.12),(2.13) and (2.30), we have
\begin{align}
\widetilde{{\rm Ric}}^1\left(\begin{array}{c}
e_1\\
e_2\\
e_3
\end{array}\right)=\left(\begin{array}{ccc}
-\left(\alpha^2+{\beta^2}\right)&\alpha\beta&\frac{\alpha\beta}{2}\\
\alpha\beta&-\left(\alpha^2+{\beta^2}\right)&-\frac{\alpha^2}{2}\\
-\frac{\alpha\beta}{2}&\frac{\alpha^2}{2}&0
\end{array}\right)\left(\begin{array}{c}
e_1\\
e_2\\
e_3
\end{array}\right).
\end{align}
If $(G_1,g,J)$ is the second kind algebraic Ricci soliton associated to the connection $\nabla^1$, then
$\widetilde{{\rm Ric}}^1=c{\rm Id}+D$, so
\begin{align}
\left\{\begin{array}{l}
De_1=-(\alpha^2+{\beta^2}+c)e_1+\alpha\beta e_2+\frac{\alpha\beta}{2}e_3,\\
De_2=\alpha\beta e_1-(\alpha^2+{\beta^2}+c)e_2-\frac{\alpha^2}{2}e_3,\\
De_3=-\frac{\alpha\beta}{2}e_1+\frac{\alpha^2}{2}e_2-ce_3.\\
\end{array}\right.
\end{align}
By (2.15) and (2.34), we get $\frac{\alpha^2}{2}+c=0$, $\beta=0$. Then
\vskip 0.5 true cm
\begin{thm}
$(G_1,g,J)$ is the second kind algebraic Ricci soliton associated to the connection $\nabla^1$ if and only if  $\frac{\alpha^2}{2}+c=0$, $\beta=0$, $\alpha\neq 0$.
In particular,
\begin{align}
\widetilde{{\rm Ric}}^1\left(\begin{array}{c}
e_1\\
e_2\\
e_3
\end{array}\right)=\left(\begin{array}{ccc}
-\alpha^2&0&0\\
0&-\alpha^2&-\frac{\alpha^2}{2}\\
0&\frac{\alpha^2}{2}&0
\end{array}\right)\left(\begin{array}{c}
e_1\\
e_2\\
e_3
\end{array}\right),~~~
{D}\left(\begin{array}{c}
e_1\\
e_2\\
e_3
\end{array}\right)=\left(\begin{array}{ccc}
-\frac{\alpha^2}{2}&0&0\\
0&-\frac{\alpha^2}{2}&-\frac{\alpha^2}{2}\\
0&\frac{\alpha^2}{2}&\frac{\alpha^2}{2}
\end{array}\right)\left(\begin{array}{c}
e_1\\
e_2\\
e_3
\end{array}\right).
\end{align}
\end{thm}

\vskip 0.5 true cm
\noindent{\bf 2.2 Algebraic Ricci solitons of $G_2$}\\
\vskip 0.5 true cm
By (2.2) and Lemma 3.5 in \cite{BO}, we have for $G_2$, there exists a pseudo-orthonormal basis $\{e_1,e_2,e_3\}$ with $e_3$ timelike such that the Lie
algebra of $G_2$ satisfies
\begin{equation}
[e_1,e_2]=\gamma e_2-\beta e_3,~~[e_1,e_3]=-\beta e_2-\gamma e_3,~~[e_2,e_3]=\alpha e_1,~~\gamma\neq 0.
\end{equation}
\vskip 0.5 true cm
\begin{lem}(\cite{Ca},\cite{BO})
The Levi-Civita connection $\nabla$ of $G_2$ is given by
\begin{align}
&\nabla_{e_1}e_1=0,~~\nabla_{e_2}e_1=-\gamma e_2+\frac{\alpha}{2}e_3,~~\nabla_{e_3}e_1=\frac{\alpha}{2}e_2+\gamma e_3,\\\notag
&\nabla_{e_1}e_2=(\frac{\alpha}{2}-\beta) e_3,~~\nabla_{e_2}e_2=\gamma e_1,~~\nabla_{e_3}e_2=-\frac{\alpha}{2}e_1,\\\notag
&\nabla_{e_1}e_3=(\frac{\alpha}{2}-\beta) e_2,~~\nabla_{e_2}e_3=\frac{\alpha}{2}e_1,~~\nabla_{e_3}e_3=\gamma e_1.
\notag
\end{align}
\end{lem}
By (2.4) and (2.37), we have
\vskip 0.5 true cm
\begin{lem}
For $G_2$, the following equalities hold
\begin{align}
&\nabla_{e_1}(J)e_1=0,~~\nabla_{e_1}(J)e_2=(\alpha-2\beta)e_3,~~\nabla_{e_1}(J)e_3=-(\alpha-2\beta)e_2,~~\\\notag
&\nabla_{e_2}(J)e_1=\alpha e_3,~~\nabla_{e_2}(J)e_2=0,~~\nabla_{e_2}(J)e_3=-\alpha e_1,~~\\\notag
&\nabla_{e_3}(J)e_1=2\gamma e_3,~~\nabla_{e_3}(J)e_2=0,~~\nabla_{e_3}(J)e_3=-2\gamma e_1.~~\notag
\end{align}
\end{lem}
\vskip 0.5 true cm
By (2.4),(2.5), Lemma 2.12 and Lemma 2.13, we have
\vskip 0.5 true cm
\begin{lem}
The canonical connection $\nabla^0$ of $(G_2,J)$ is given by
\begin{align}
&\nabla^0_{e_1}e_1=0,~~\nabla^0_{e_1}e_2=0,~~\nabla^0_{e_1}e_3=0,\\\notag
&\nabla^0_{e_2}e_1=-\gamma e_2,~~\nabla^0_{e_2}e_2=\gamma e_1,~~\nabla^0_{e_2}e_3=0,\\\notag
&\nabla^0_{e_3}e_1=\frac{\alpha}{2}e_2,~~\nabla^0_{e_3}e_2=-\frac{\alpha}{2}e_1,~~\nabla^0_{e_3}e_3=0.
\notag
\end{align}
\end{lem}
\vskip 0.5 true cm
\indent By (2.7) and Lemma 2.14, we have
\vskip 0.5 true cm
\begin{lem}
The curvature $R^0$ of the canonical connection $\nabla^0$ of $(G_2,J)$ is given by
\begin{align}
&R^0(e_1,e_2)e_1=(\gamma^2+\frac{\alpha\beta}{2})e_2,~~R^0(e_1,e_2)e_2=-(\gamma^2+\frac{\alpha\beta}{2})e_1,~~R^0(e_1,e_2)e_3=0,\\\notag
&R^0(e_1,e_3)e_1=(\frac{\alpha\gamma}{2}-\beta\gamma) e_2,~~R^0(e_1,e_3)e_2=-(\frac{\alpha\gamma}{2}-\beta\gamma) e_1,~~R^0(e_1,e_3)e_3=0,\\\notag
&R^0(e_2,e_3)e_1=0,~~R^0(e_2,e_3)e_2=0,~~R^0(e_2,e_3)e_3=0.\notag
\notag
\end{align}
\end{lem}
\vskip 0.5 true cm
By (2.9),(2.11) and Lemma 2.15, we get for $(G_2,\nabla^0)$
\begin{align}
{\rm Ric}^0\left(\begin{array}{c}
e_1\\
e_2\\
e_3
\end{array}\right)=\left(\begin{array}{ccc}
-\left(\gamma^2+\frac{\alpha\beta}{2}\right)&0&0\\
0&-\left(\gamma^2+\frac{\alpha\beta}{2}\right)&0\\
0&\beta\gamma-\frac{\alpha\gamma}{2}&0
\end{array}\right)\left(\begin{array}{c}
e_1\\
e_2\\
e_3
\end{array}\right).
\end{align}
If $(G_2,g,J)$ is the first kind algebraic Ricci soliton associated to the connection $\nabla^0$, then
${\rm Ric}^0=c{\rm Id}+D$, so
\begin{align}
\left\{\begin{array}{l}
De_1=-(\gamma^2+\frac{\alpha\beta}{2}+c)e_1,\\
De_2=-(\gamma^2+\frac{\alpha\beta}{2}+c)e_2,\\
De_3=(\beta\gamma-\frac{\alpha\gamma}{2})e_2-ce_3.\\
\end{array}\right.
\end{align}
By (2.15) and (2.42), we get $\alpha=\beta=0$. $\gamma^2+c=0$. Then
\vskip 0.5 true cm
\begin{thm}
$(G_2,g,J)$ is the first kind algebraic Ricci soliton associated to the connection $\nabla^0$ if and only if $\alpha=\beta=0$. $\gamma^2+c=0$, $\gamma\neq 0$.
In particular,
\begin{align}
{\rm Ric}^0\left(\begin{array}{c}
e_1\\
e_2\\
e_3
\end{array}\right)=\left(\begin{array}{ccc}
-\gamma^2&0&0\\
0&-\gamma^2&0\\
0&0&0
\end{array}\right)\left(\begin{array}{c}
e_1\\
e_2\\
e_3
\end{array}\right),~~~
{D}\left(\begin{array}{c}
e_1\\
e_2\\
e_3
\end{array}\right)=\left(\begin{array}{ccc}
0&0&0\\
0&0&0\\
0&0&\gamma^2
\end{array}\right)\left(\begin{array}{c}
e_1\\
e_2\\
e_3
\end{array}\right).
\end{align}
\end{thm}
\vskip 0.5 true cm
\indent By (2.12),(2.13) and (2.41), we have for $(G_2,\nabla^0)$
\begin{align}
\widetilde{{\rm Ric}}^0\left(\begin{array}{c}
e_1\\
e_2\\
e_3
\end{array}\right)=\left(\begin{array}{ccc}
-\left(\gamma^2+\frac{\alpha\beta}{2}\right)&0&0\\
0&-\left(\gamma^2+\frac{\alpha\beta}{2}\right)&\frac{\alpha\gamma}{4}-\frac{\beta\gamma}{2}\\
0&\frac{\beta\gamma}{2}-\frac{\alpha\gamma}{4}&0
\end{array}\right)\left(\begin{array}{c}
e_1\\
e_2\\
e_3
\end{array}\right).
\end{align}
If $(G_2,g,J)$ is the second kind algebraic Ricci soliton associated to the connection $\nabla^0$, then
$\widetilde{{\rm Ric}}^0=c{\rm Id}+D$, so
\begin{align}
\left\{\begin{array}{l}
De_1=-(\gamma^2+\frac{\alpha\beta}{2}+c)e_1,\\
De_2=-(\gamma^2+\frac{\alpha\beta}{2}+c)e_2+(\frac{\alpha\gamma}{4}-\frac{\beta\gamma}{2})e_3,\\
De_3=(\frac{\beta\gamma}{2}-\frac{\alpha\gamma}{4})e_2-ce_3.\\
\end{array}\right.
\end{align}
By (2.15) and (2.45), we get
\begin{align}
\left\{\begin{array}{l}
\gamma^2+\alpha\beta-\beta^2+c=0,\\
\beta(2\gamma^2+\alpha\beta+c)+\gamma(\beta\gamma-\frac{\alpha\gamma}{2})=0,\\
\beta c+\gamma(\beta\gamma-\frac{\alpha\gamma}{2})=0,\\
\alpha c=0.\\
\end{array}\right.
\end{align}
Solving (2.46), we get $\alpha=\beta=0$, $\gamma^2+c=0$.
Then
\vskip 0.5 true cm
\begin{thm}
$(G_2,g,J)$ is the second kind algebraic Ricci soliton associated to the connection $\nabla^0$ if and only if  $\alpha=\beta=0$, $\gamma^2+c=0$, $\gamma\neq 0$.
In particular,
\begin{align}
\widetilde{{\rm Ric}}^0\left(\begin{array}{c}
e_1\\
e_2\\
e_3
\end{array}\right)=\left(\begin{array}{ccc}
-\gamma^2&0&0\\
0&-\gamma^2&0\\
0&0&0
\end{array}\right)\left(\begin{array}{c}
e_1\\
e_2\\
e_3
\end{array}\right),~~~
{D}\left(\begin{array}{c}
e_1\\
e_2\\
e_3
\end{array}\right)=\left(\begin{array}{ccc}
0&0&0\\
0&0&0\\
0&0&\gamma^2
\end{array}\right)\left(\begin{array}{c}
e_1\\
e_2\\
e_3
\end{array}\right).
\end{align}
\end{thm}
\vskip 0.5 true cm
\indent By (2.6), Lemma 2.13 and Lemma 2.14, we have
\vskip 0.5 true cm
\begin{lem}
The Kobayashi-Nomizu connection $\nabla^1$ of $(G_2,J)$ is given by
\begin{align}
&\nabla^1_{e_1}e_1=0,~~\nabla^1_{e_1}e_2=0,~~\nabla^1_{e_1}e_3=-\gamma e_3,\\\notag
&\nabla^1_{e_2}e_1=-\gamma e_2,~~\nabla^1_{e_2}e_2=\gamma e_1,~~\nabla^1_{e_2}e_3=0,\\\notag
&\nabla^1_{e_3}e_1=\beta e_2,~~\nabla^1_{e_3}e_2=-\alpha e_1,~~\nabla^1_{e_3}e_3=0.
\notag
\end{align}
\end{lem}
\indent By (2.8) and Lemma 2.18, we have
\vskip 0.5 true cm
\begin{lem}
The curvature $R^1$ of the Kobayashi-Nomizu connection $\nabla^1$ of $(G_2,J)$ is given by
\begin{align}
&R^1(e_1,e_2)e_1=(\beta^2+{\gamma^2})e_2,~~R^1(e_1,e_2)e_2=-(\gamma^2+\alpha\beta)e_1,~~R^1(e_1,e_2)e_3=0,\\\notag
&R^1(e_1,e_3)e_1=0,~~R^1(e_1,e_3)e_2=(\beta\gamma-\alpha\gamma) e_1,~~R^1(e_1,e_3)e_3=0,\\\notag
&R^1(e_2,e_3)e_1=(\beta\gamma-\alpha\gamma) e_1,~~R^1(e_2,e_3)e_2=(\alpha\gamma-\beta\gamma) e_2,~~R^1(e_2,e_3)e_3=\alpha\gamma e_3.\notag
\notag
\end{align}
\end{lem}
\vskip 0.5 true cm
By (2.10),(2.11) and Lemma 2.19, we get for $(G_2,\nabla^1)$
\begin{align}
{\rm Ric}^1\left(\begin{array}{c}
e_1\\
e_2\\
e_3
\end{array}\right)=\left(\begin{array}{ccc}
-\left(\beta^2+\gamma^2\right)&0&0\\
0&-(\gamma^2+\alpha\beta)&\alpha\gamma\\
0&0&0
\end{array}\right)\left(\begin{array}{c}
e_1\\
e_2\\
e_3
\end{array}\right).
\end{align}
If $(G_2,g,J)$ is the first kind algebraic Ricci soliton associated to the connection $\nabla^1$, then
${\rm Ric}^1=c{\rm Id}+D$, so
\begin{align}
\left\{\begin{array}{l}
De_1=-\left(\beta^2+\gamma^2+c\right)e_1,\\
De_2=-(\gamma^2+\alpha\beta+c)e_2+\alpha\gamma e_3,\\
De_3=-ce_3.\\
\end{array}\right.
\end{align}
By (2.15) and (2.51), we get
\begin{align}
\left\{\begin{array}{l}
\beta^2+\gamma^2+c+\alpha\beta=0,\\
\beta(\beta^2+2\gamma^2+c+\alpha\beta)-2\alpha\gamma^2=0,\\
\beta(\beta^2-\alpha\beta+c)=0,\\
\alpha(\beta^2-\alpha\beta-c)=0.\\
\end{array}\right.
\end{align}
Solving (2.52), we get $\alpha=\beta=0$, $\gamma^2+c=0$.
Then
\vskip 0.5 true cm
\begin{thm}
$(G_2,g,J)$ is the first kind algebraic Ricci soliton associated to the connection $\nabla^1$ if and only if $\alpha=\beta=0$, $\gamma^2+c=0$, $\gamma\neq 0$.
In particular,
\begin{align}
{\rm Ric}^1\left(\begin{array}{c}
e_1\\
e_2\\
e_3
\end{array}\right)=\left(\begin{array}{ccc}
-\gamma^2&0&0\\
0&-\gamma^2&0\\
0&0&0
\end{array}\right)\left(\begin{array}{c}
e_1\\
e_2\\
e_3
\end{array}\right),~~~
{D}\left(\begin{array}{c}
e_1\\
e_2\\
e_3
\end{array}\right)=\left(\begin{array}{ccc}
0&0&0\\
0&0&0\\
0&0&\gamma^2
\end{array}\right)\left(\begin{array}{c}
e_1\\
e_2\\
e_3
\end{array}\right).
\end{align}
\end{thm}
\vskip 0.5 true cm
\indent By (2.12),(2.13) and (2.50), we have for $(G_2,\nabla^1)$
\begin{align}
\widetilde{{\rm Ric}}^1\left(\begin{array}{c}
e_1\\
e_2\\
e_3
\end{array}\right)=\left(\begin{array}{ccc}
-\left(\beta^2+\gamma^2\right)&0&0\\
0&-\left(\gamma^2+\alpha\beta\right)&\frac{\alpha\gamma}{2}\\
0&-\frac{\alpha\gamma}{2}&0
\end{array}\right)\left(\begin{array}{c}
e_1\\
e_2\\
e_3
\end{array}\right).
\end{align}
If $(G_2,g,J)$ is the second kind algebraic Ricci soliton associated to the connection $\nabla^1$, then
$\widetilde{{\rm Ric}}^1=c{\rm Id}+D$, so
\begin{align}
\left\{\begin{array}{l}
De_1=-(\beta^2+\gamma^2+c)e_1,\\
De_2=-(\gamma^2+\alpha\beta+c)e_2+\frac{\alpha\gamma}{2}e_3,\\
De_3=-\frac{\alpha\gamma}{2}e_2-ce_3.\\
\end{array}\right.
\end{align}
By (2.15) and (2.55), we get
\begin{align}
\left\{\begin{array}{l}
\beta^2+\gamma^2+c+\alpha\beta=0,\\
\beta(\beta^2+2\gamma^2+c+\alpha\beta)-\alpha\gamma^2=0,\\
\beta(\beta^2-\alpha\beta+c)-\alpha\gamma^2=0,\\
\alpha(\beta^2-\alpha\beta-c)=0.\\
\end{array}\right.
\end{align}
Solving (2.56), we get $\alpha=\beta=0$, $\gamma^2+c=0$. Then
\vskip 0.5 true cm
\begin{thm}
$(G_2,g,J)$ is the second kind algebraic Ricci soliton associated to the connection $\nabla^1$ if and only if $\alpha=\beta=0$, $\gamma^2+c=0$,  $\gamma\neq 0$.
In particular,
\begin{align}
\widetilde{{\rm Ric}}^1\left(\begin{array}{c}
e_1\\
e_2\\
e_3
\end{array}\right)=\left(\begin{array}{ccc}
-\gamma^2&0&0\\
0&-\gamma^2&0\\
0&0&0
\end{array}\right)\left(\begin{array}{c}
e_1\\
e_2\\
e_3
\end{array}\right),~~~
{D}\left(\begin{array}{c}
e_1\\
e_2\\
e_3
\end{array}\right)=\left(\begin{array}{ccc}
0&0&0\\
0&0&0\\
0&0&\gamma^2
\end{array}\right)\left(\begin{array}{c}
e_1\\
e_2\\
e_3
\end{array}\right).
\end{align}
\end{thm}
\vskip 0.5 true cm
\noindent{\bf 2.3 Algebraic Ricci solitons of $G_3$}\\
\vskip 0.5 true cm
By (2.3) and Lemma 3.8 in \cite{BO}, we have for $G_3$, there exists a pseudo-orthonormal basis $\{e_1,e_2,e_3\}$ with $e_3$ timelike such that the Lie
algebra of $G_3$ satisfies
\begin{equation}
[e_1,e_2]=-\gamma e_3,~~[e_1,e_3]=-\beta e_2,~~[e_2,e_3]=\alpha e_1.
\end{equation}
\vskip 0.5 true cm
\begin{lem}(\cite{Ca},\cite{BO})
The Levi-Civita connection $\nabla$ of $G_3$ is given by
\begin{align}
&\nabla_{e_1}e_1=0,~~\nabla_{e_2}e_1=a_2 e_3,~~\nabla_{e_3}e_1=a_3e_2,\\\notag
&\nabla_{e_1}e_2=a_1 e_3,~~\nabla_{e_2}e_2=0,~~\nabla_{e_3}e_2=-a_3e_1,\\\notag
&\nabla_{e_1}e_3=a_1 e_2,~~\nabla_{e_2}e_3=a_2e_1,~~\nabla_{e_3}e_3=0,
\notag
\end{align}
where
\begin{equation}
a_1=\frac{1}{2}(\alpha-\beta-\gamma),~~a_2=\frac{1}{2}(\alpha-\beta+\gamma),~~a_3=\frac{1}{2}(\alpha+\beta-\gamma).
\end{equation}
\end{lem}
By (2.4) and (2.59), we have
\vskip 0.5 true cm
\begin{lem}
For $G_3$, the following equalities hold
\begin{align}
&\nabla_{e_1}(J)e_1=0,~~\nabla_{e_1}(J)e_2=2a_1e_3,~~\nabla_{e_1}(J)e_3=-2a_1e_2,~~\\\notag
&\nabla_{e_2}(J)e_1=2a_2 e_3,~~\nabla_{e_2}(J)e_2=0,~~\nabla_{e_2}(J)e_3=-2a_2 e_1,~~\\\notag
&\nabla_{e_3}(J)e_1=0,~~\nabla_{e_3}(J)e_2=0,~~\nabla_{e_3}(J)e_3=0.~~\notag
\end{align}
\end{lem}
\vskip 0.5 true cm
By (2.4),(2.5), Lemma 2.22 and Lemma 2.23, we have
\vskip 0.5 true cm
\begin{lem}
The canonical connection $\nabla^0$ of $(G_3,J)$ is given by
\begin{align}
&\nabla^0_{e_1}e_1=0,~~\nabla^0_{e_1}e_2=0,~~\nabla^0_{e_1}e_3=0,\\\notag
&\nabla^0_{e_2}e_1=0,~~\nabla^0_{e_2}e_2=0,~~\nabla^0_{e_2}e_3=0,\\\notag
&\nabla^0_{e_3}e_1=a_3e_2,~~\nabla^0_{e_3}e_2=-a_3e_1,~~\nabla^0_{e_3}e_3=0.
\notag
\end{align}
\end{lem}
\vskip 0.5 true cm
\indent By (2.7) and Lemma 2.24, we have
\vskip 0.5 true cm
\begin{lem}
The curvature $R^0$ of the canonical connection $\nabla^0$ of $(G_3,J)$ is given by
\begin{align}
&R^0(e_1,e_2)e_1=\gamma a_3e_2,~~R^0(e_1,e_2)e_2=-\gamma a_3e_1,~~R^0(e_1,e_2)e_3=0,\\\notag
&R^0(e_1,e_3)e_1=0,~~R^0(e_1,e_3)e_2=0,~~R^0(e_1,e_3)e_3=0,\\\notag
&R^0(e_2,e_3)e_1=0,~~R^0(e_2,e_3)e_2=0,~~R^0(e_2,e_3)e_3=0.\notag
\notag
\end{align}
\end{lem}
\vskip 0.5 true cm
By (2.9),(2.11) and Lemma 2.25, we get for $(G_3,\nabla^0)$
\begin{align}
{\rm Ric}^0\left(\begin{array}{c}
e_1\\
e_2\\
e_3
\end{array}\right)=\left(\begin{array}{ccc}
-\gamma a_3&0&0\\
0&-\gamma a_3&0\\
0&0&0
\end{array}\right)\left(\begin{array}{c}
e_1\\
e_2\\
e_3
\end{array}\right).
\end{align}
If $(G_3,g,J)$ is the first kind algebraic Ricci soliton associated to the connection $\nabla^0$, then
${\rm Ric}^0=c{\rm Id}+D$, so
\begin{align}
\left\{\begin{array}{l}
De_1=-(\gamma a_3+c)e_1,\\
De_2=-(\gamma a_3+c)e_2,\\
De_3=-ce_3.\\
\end{array}\right.
\end{align}
By (2.15) and (2.65), we get
\begin{align}
\left\{\begin{array}{l}
\gamma(2\gamma a_3+c)=0,\\
\beta c=0,\\
\alpha c=0.\\
\end{array}\right.
\end{align}

Solving (2.66), then we have
\vskip 0.5 true cm
\begin{thm}
$(G_3,g,J)$ is the first kind algebraic Ricci soliton associated to the connection $\nabla^0$ if and only if \\
(i)$\alpha=\beta=\gamma=0$, in particular
\begin{align}
{D}\left(\begin{array}{c}
e_1\\
e_2\\
e_3
\end{array}\right)=\left(\begin{array}{ccc}
-c&0&0\\
0&-c&0\\
0&0&-c
\end{array}\right)\left(\begin{array}{c}
e_1\\
e_2\\
e_3 \notag \\
\end{array}\right).
\end{align}
(ii)$\alpha=\beta=0$, $\gamma^2=c$, in particular
\begin{align}
{D}\left(\begin{array}{c}
e_1\\
e_2\\
e_3 \notag \\
\end{array}\right)=\left(\begin{array}{ccc}
-\frac{\gamma^2}{2}&0&0\\
0&-\frac{\gamma^2}{2}&0\\
0&0&-\gamma^2
\end{array}\right)\left(\begin{array}{c}
e_1\\
e_2\\
e_3 \notag \\
\end{array}\right).
\end{align}
(iii)$\alpha\neq 0$ or $\beta\neq 0$, $\gamma=0$, $c=0$, in particular
\begin{align}
{D}\left(\begin{array}{c}
e_1\\
e_2\\
e_3
\end{array}\right)=\left(\begin{array}{ccc}
0&0&0\\
0&0&0\\
0&0&0
\end{array}\right)\left(\begin{array}{c}
e_1\\
e_2\\
e_3 \notag \\
\end{array}\right).
\end{align}
(iv)$\alpha\neq 0$ or $\beta\neq 0$, $\gamma=\alpha+\beta$, $c=0$, in particular
\begin{align}
{D}\left(\begin{array}{c}
e_1\\
e_2\\
e_3
\end{array}\right)=\left(\begin{array}{ccc}
0&0&0\\
0&0&0\\
0&0&0
\end{array}\right)\left(\begin{array}{c}
e_1\\
e_2\\
e_3  \notag \\
\end{array}\right).
\end{align}
\end{thm}
\vskip 0.5 true cm
Since $\rho^0(e_i,e_j)=\rho^0(e_j,e_i)$, then $\widetilde{\rho}^0(e_i,e_j)=\rho^0(e_i,e_j)$. Then
$(G_3,g,J)$ is the second kind algebraic Ricci soliton associated to the connection $\nabla^0$ if and only if
$(G_3,g,J)$ is the first kind algebraic Ricci soliton associated to the connection $\nabla^0$.
By (2.4), (2.6), Lemma 2.23 and Lemma 2.24, we have
\begin{lem}
The Kobayashi-Nomizu connection $\nabla^1$ of $(G_3,J)$ is given by
\begin{align}
&\nabla^1_{e_1}e_1=0,~~\nabla^1_{e_1}e_2=0,~~\nabla^1_{e_1}e_3=0,\\\notag
&\nabla^1_{e_2}e_1=0,~~\nabla^1_{e_2}e_2=0,~~\nabla^1_{e_2}e_3=0,\\\notag
&\nabla^1_{e_3}e_1=(a_3-a_1)e_2,~~\nabla^1_{e_3}e_2=-(a_2+a_3) e_1,~~\nabla^1_{e_3}e_3=0.
\notag
\end{align}
\end{lem}
\indent By (2.8) and Lemma 2.27, we have
\vskip 0.5 true cm
\begin{lem}
The curvature $R^1$ of the Kobayashi-Nomizu connection $\nabla^1$ of $(G_3,J)$ is given by
\begin{align}
&R^1(e_1,e_2)e_1=\gamma (a_3-a_1)e_2,~~R^1(e_1,e_2)e_2=-\gamma (a_2+a_3)e_1,~~R^1(e_1,e_2)e_3=0,\\\notag
&R^1(e_1,e_3)e_1=0,~~R^1(e_1,e_3)e_2=0,~~R^1(e_1,e_3)e_3=0,\\\notag
&R^1(e_2,e_3)e_1=0,~~R^1(e_2,e_3)e_2=0,~~R^1(e_2,e_3)e_3=0.\notag
\notag
\end{align}
\end{lem}
\vskip 0.5 true cm
By (2.10),(2.11) and Lemma 2.28, we get for $(G_3,\nabla^1)$
\begin{align}
{\rm Ric}^1\left(\begin{array}{c}
e_1\\
e_2\\
e_3
\end{array}\right)=\left(\begin{array}{ccc}
\gamma(a_1-a_3)&0&0\\
0&-\gamma(a_2+a_3)&0\\
0&0&0
\end{array}\right)\left(\begin{array}{c}
e_1\\
e_2\\
e_3
\end{array}\right).
\end{align}
If $(G_3,g,J)$ is the first kind algebraic Ricci soliton associated to the connection $\nabla^1$, then
${\rm Ric}^1=c{\rm Id}+D$, so
\begin{align}
\left\{\begin{array}{l}
De_1=\left[\gamma(a_1-a_3)-c\right]e_1,\\
De_2=[-\gamma (a_2+a_3)-c]e_2,\\
De_3=-ce_3.\\
\end{array}\right.
\end{align}
By (2.15) and (2.70), we get
\begin{align}
\left\{\begin{array}{l}
\gamma[\gamma(\alpha+\beta]+c]=0,\\
\beta[\gamma(\alpha-\beta)-c]=0,\\
\alpha[\gamma(\alpha-\beta)+c]=0,
\end{array}\right.
\end{align}
Solving (2.71), we get
\vskip 0.5 true cm
\begin{thm}
$(G_3,g,J)$ is the first kind algebraic Ricci soliton associated to the connection $\nabla^1$ if and only if \\
(i) $\alpha\beta\neq 0$, $\gamma=0$, $c=0$,
in particular,
\begin{align}
{D}\left(\begin{array}{c}
e_1\\
e_2\\
e_3
\end{array}\right)=\left(\begin{array}{ccc}
0&0&0\\
0&0&0\\
0&0&0
\end{array}\right)\left(\begin{array}{c}
e_1\\
e_2\\
e_3 \notag \\
\end{array}\right).
\end{align}
(ii)$\alpha=\beta=\gamma=0$, $c\neq 0$,
in particular,
 \begin{align}
{D}\left(\begin{array}{c}
e_1\\
e_2\\
e_3
\end{array}\right)=\left(\begin{array}{ccc}
-c&0&0\\
0&-c&0\\
0&0&-c
\end{array}\right)\left(\begin{array}{c}
e_1\\
e_2\\
e_3 \notag \\
\end{array}\right).
\end{align}
(iii)$\alpha=0$, $\gamma\beta+c=0$,
in particular,
\begin{align}
{D}\left(\begin{array}{c}
e_1\\
e_2\\
e_3
\end{array}\right)=\left(\begin{array}{ccc}
0&0&0\\
0&-c&0\\
0&0&-c
\end{array}\right)\left(\begin{array}{c}
e_1\\
e_2\\
e_3 \notag \\
\end{array}\right).
\end{align}
(iv)$\beta=0$, $\gamma\alpha+c=0$,
in particular,
 \begin{align}
{D}\left(\begin{array}{c}
e_1\\
e_2\\
e_3
\end{array}\right)=\left(\begin{array}{ccc}
-c&0&0\\
0&0&0\\
0&0&-c
\end{array}\right)\left(\begin{array}{c}
e_1\\
e_2\\
e_3 \notag \\
\end{array}\right).
\end{align}
\end{thm}
\vskip 0.5 true cm
Since $\rho^1(e_i,e_j)=\rho^1(e_j,e_i)$, then $\widetilde{\rho}^1(e_i,e_j)=\rho^1(e_i,e_j)$. Then
$(G_3,g,J)$ is the second kind algebraic Ricci soliton associated to the connection $\nabla^1$ if and only if
$(G_3,g,J)$ is the first kind algebraic Ricci soliton associated to the connection $\nabla^1$.
\vskip 0.5 true cm
\noindent{\bf 2.4 Algebraic Ricci solitons of $G_4$}
\vskip 0.5 true cm
By (2.4) and Lemma 3.11 in \cite{BO}, we have for $G_4$, there exists a pseudo-orthonormal basis $\{e_1,e_2,e_3\}$ with $e_3$ timelike such that the Lie
algebra of $G_4$ satisfies
\begin{align}
[e_1,e_2]=-e_2+(2\eta-\beta)e_3,~~\eta=1~{\rm or}-1,~~[e_1,e_3]=-\beta e_2+ e_3,~~[e_2,e_3]=\alpha e_1.
\end{align}
\vskip 0.5 true cm
\begin{lem}(\cite{Ca},\cite{BO})
The Levi-Civita connection $\nabla$ of $G_4$ is given by
\begin{align}
&\nabla_{e_1}e_1=0,~~\nabla_{e_2}e_1=e_2+b_2e_3,~~\nabla_{e_3}e_1=b_3e_2-e_3,\\\notag
&\nabla_{e_1}e_2=b_1 e_3,~~\nabla_{e_2}e_2=- e_1,~~\nabla_{e_3}e_2=-b_3e_1,\\\notag
&\nabla_{e_1}e_3=b_1 e_2,~~\nabla_{e_2}e_3=b_2e_1,~~\nabla_{e_3}e_3=- e_1,
\notag
\end{align}
where
\begin{equation}
b_1=\frac{\alpha}{2}+\eta-\beta,~~b_2=\frac{\alpha}{2}-\eta,~~b_3=\frac{\alpha}{2}+\eta.
\end{equation}
\end{lem}
\vskip 0.5 true cm
By (2.4) and Lemma 2.30, we have
\vskip 0.5 true cm
\begin{lem}
For $G_4$, the following equalities hold
\begin{align}
&\nabla_{e_1}(J)e_1=0,~~\nabla_{e_1}(J)e_2=2b_1e_3,~~\nabla_{e_1}(J)e_3=-2b_1e_2,~~\\\notag
&\nabla_{e_2}(J)e_1=2b_2 e_3,~~\nabla_{e_2}(J)e_2=0,~~\nabla_{e_2}(J)e_3=-2b_2 e_1,~~\\\notag
&\nabla_{e_3}(J)e_1=-2e_3,~~\nabla_{e_3}(J)e_2=0,~~\nabla_{e_3}(J)e_3=2 e_1.~~\notag
\end{align}
\end{lem}
\vskip 0.5 true cm
By (2.4),(2.5), Lemma 2.30 and Lemma 2.31, we have
\vskip 0.5 true cm
\begin{lem}
The canonical connection $\nabla^0$ of $(G_4,J)$ is given by
\begin{align}
&\nabla^0_{e_1}e_1=0,~~\nabla^0_{e_1}e_2=0,~~\nabla^0_{e_1}e_3=0,\\\notag
&\nabla^0_{e_2}e_1= e_2,~~\nabla^0_{e_2}e_2=- e_1,~~\nabla^0_{e_2}e_3=0,\\\notag
&\nabla^0_{e_3}e_1=b_3e_2,~~\nabla^0_{e_3}e_2=-b_3e_1,~~\nabla^0_{e_3}e_3=0.
\notag
\end{align}
\end{lem}
\vskip 0.5 true cm
\indent By (2.7) and Lemma 2.32, we have
\vskip 0.5 true cm
\begin{lem}
The curvature $R^0$ of the canonical connection $\nabla^0$ of $(G_4,J)$ is given by
\begin{align}
&R^0(e_1,e_2)e_1=[(\beta-2\eta)b_3+1]e_2,~~R^0(e_1,e_2)e_2=[b_3(2\eta-\beta)-1]e_1,~~R^0(e_1,e_2)e_3=0,\\\notag
&R^0(e_1,e_3)e_1=(\beta-b_3) e_2,~~R^0(e_1,e_3)e_2=(b_3-\beta) e_1,~~R^0(e_1,e_3)e_3=0,\\\notag
&R^0(e_2,e_3)e_1=0,~~R^0(e_2,e_3)e_2=0,~~R^0(e_2,e_3)e_3=0.\notag
\notag
\end{align}
\end{lem}
\vskip 0.5 true cm
By (2.9),(2.11) and Lemma 2.33, we get for $(G_4,\nabla^0)$
\begin{align}
{\rm Ric}^0\left(\begin{array}{c}
e_1\\
e_2\\
e_3
\end{array}\right)=\left(\begin{array}{ccc}
[b_3(2\eta-\beta)-1]&0&0\\
0&[b_3(2\eta-\beta)-1]&0\\
0&b_3-\beta&0
\end{array}\right)\left(\begin{array}{c}
e_1\\
e_2\\
e_3
\end{array}\right).
\end{align}
If $(G_4,g,J)$ is the first kind algebraic Ricci soliton associated to the connection $\nabla^0$, then
${\rm Ric}^0=c{\rm Id}+D$, so
\begin{align}
\left\{\begin{array}{l}
De_1=[b_3(2\eta-\beta)-1-c]e_1,\\
De_2=[b_3(2\eta-\beta)-1-c]e_2,\\
De_3=(b_3-\beta)e_2-ce_3.\\
\end{array}\right.
\end{align}
By (2.15) and (2.79), we get
\begin{align}
\left\{\begin{array}{l}
(2\eta-\beta)(2b_3-\beta)-c-1=0,\\
\left(2b_3(2\eta-\beta)-2-c\right)(2\eta-\beta)=0,\\
2(b_3-\beta)-c\beta=0,\\
c\alpha=0.\
\end{array}\right.
\end{align}
Solving (2.80), we get
\vskip 0.5 true cm
\begin{thm}
$(G_4,g,J)$ is the first kind algebraic Ricci soliton associated to the connection $\nabla^0$ if and only if \\
(i) $\alpha=0$, $\beta=1$, $c=0$, $\eta=1$,
in particular,
\begin{align}
{D}\left(\begin{array}{c}
e_1\\
e_2\\
e_3 \notag \\
\end{array}\right)=\left(\begin{array}{ccc}
0&0&0\\
0&0&0\\
0&0&0
\end{array}\right)\left(\begin{array}{c}
e_1\\
e_2\\
e_3
\end{array}\right).
\end{align}
(ii)$\alpha=0$, $c=-1$, $\beta=2\eta$,
 \begin{align}
{D}\left(\begin{array}{c}
e_1\\
e_2\\
e_3
\end{array}\right)=\left(\begin{array}{ccc}
0&0&0\\
0&0&0\\
0&-\eta&1
\end{array}\right)\left(\begin{array}{c}
e_1\\
e_2\\
e_3 \notag \\
\end{array}\right).
\end{align}
\end{thm}
\vskip 0.5 true cm
\indent By (2.12),(2.13) and (2.78), we have for $(G_4,\nabla^0)$
\begin{align}
\widetilde{{\rm Ric}}^0\left(\begin{array}{c}
e_1\\
e_2\\
e_3
\end{array}\right)=\left(\begin{array}{ccc}
[(2\eta-\beta)b_3-1]&0&0\\
0&[(2\eta-\beta)b_3-1]&\frac{\beta-b_3}{2}\\
0&\frac{b_3-\beta}{2}&0
\end{array}\right)\left(\begin{array}{c}
e_1\\
e_2\\
e_3
\end{array}\right).
\end{align}
If $(G_4,g,J)$ is the second kind algebraic Ricci soliton associated to the connection $\nabla^0$, then
$\widetilde{{\rm Ric}}^0=c{\rm Id}+D$, so
\begin{align}
\left\{\begin{array}{l}
De_1=[(2\eta-\beta)b_3-1-c]e_1,\\
De_2=[(2\eta-\beta)b_3-1-c]e_2+\frac{\beta-b_3}{2}e_3,\\
De_3=\frac{b_3-\beta}{2}e_2-ce_3.\\
\end{array}\right.
\end{align}
By (2.15) and (2.82), we get
\begin{align}
\left\{\begin{array}{l}
(2\eta-\beta)b_3-1-c+(\eta-\beta)(b_3-\beta)=0,\\
\beta-b_3+[2(2\eta-\beta)b_3-2-c](2\eta-\beta)=0,\\
b_3-\beta=c\beta,\\
\alpha c=0.\\
\end{array}\right.
\end{align}
Solving (2.83), we get
\vskip 0.5 true cm
\begin{thm}
$(G_4,g,J)$ is the second kind algebraic Ricci soliton associated to the connection $\nabla^0$ if and only if  $\alpha=0$, $\beta=\eta$, $c=0$.
In particular,
\begin{align}
{D}\left(\begin{array}{c}
e_1\\
e_2\\
e_3
\end{array}\right)=\left(\begin{array}{ccc}
0&0&0\\
0&0&0\\
0&0&0
\end{array}\right)\left(\begin{array}{c}
e_1\\
e_2\\
e_3 \notag \\
\end{array}\right).
\end{align}
\end{thm}
\vskip 0.5 true cm
\indent By (2.4),(2.6), Lemma 2.31 and Lemma 2.32, we have
\vskip 0.5 true cm
\begin{lem}
The Kobayashi-Nomizu connection $\nabla^1$ of $(G_4,J)$ is given by
\begin{align}
&\nabla^1_{e_1}e_1=0,~~\nabla^1_{e_1}e_2=0,~~\nabla^1_{e_1}e_3= e_3,\\\notag
&\nabla^1_{e_2}e_1= e_2,~~\nabla^1_{e_2}e_2=- e_1,~~\nabla^1_{e_2}e_3=0,\\\notag
&\nabla^0_{e_3}e_1=(b_3-b_1) e_2,~~\nabla^1_{e_3}e_2=-(b_2+b_3) e_1,~~\nabla^1_{e_3}e_3=0.
\notag
\end{align}
\end{lem}
\vskip 0.5 true cm
\indent By (2.8) and Lemma 2.36, we have
\vskip 0.5 true cm
\begin{lem}
The curvature $R^1$ of the Kobayashi-Nomizu connection $\nabla^1$ of $(G_4,J)$ is given by
\begin{align}
&R^1(e_1,e_2)e_1=[1+(\beta-2\eta)(b_3-b_1)]e_2,~~R^1(e_1,e_2)e_2=-[1+(\beta-2\eta)(b_2+b_3)]e_1,\\\notag
&R^1(e_1,e_2)e_3=0,~
R^1(e_1,e_3)e_1=(\beta-b_3+b_1)e_2,~~R^1(e_1,e_3)e_2=(b_2+b_3-\beta) e_1,\\\notag
&R^1(e_1,e_3)e_3=0,~~
R^1(e_2,e_3)e_1=(b_1+b_2) e_1,\\\notag
&R^1(e_2,e_3)e_2=-(b_1+b_2) e_2,~~R^1(e_2,e_3)e_3=-\alpha e_3.\notag
\notag
\end{align}
\end{lem}
\vskip 0.5 true cm
By (2.10),(2.11) and Lemma 2.37, we get for $(G_4,\nabla^1)$
\begin{align}
{\rm Ric}^1\left(\begin{array}{c}
e_1\\
e_2\\
e_3
\end{array}\right)=\left(\begin{array}{ccc}
-[1+(\beta-2\eta)(b_3-b_1)]&0&0\\
0&-[1+(\beta-2\eta)(b_2+b_3)]&-\alpha\\
0&(b_3-b_1-\beta)&0
\end{array}\right)\left(\begin{array}{c}
e_1\\
e_2\\
e_3
\end{array}\right).
\end{align}
If $(G_4,g,J)$ is the first kind algebraic Ricci soliton associated to the connection $\nabla^1$, then
${\rm Ric}^1=c{\rm Id}+D$, so
\begin{align}
\left\{\begin{array}{l}
De_1=-[1+(\beta-2\eta)(b_3-b_1)+c]e_1,\\
De_2=-[1+(\beta-2\eta)(b_2+b_3)+c]e_2-\alpha e_3,\\
De_3=(b_3-b_1-\beta)e_2-ce_3.\\
\end{array}\right.
\end{align}
By (2.15) and (2.87), we get
\begin{align}
\left\{\begin{array}{l}
(2\eta-\beta)(2b_3-2b_1-\beta)-c-1-\alpha\beta=0,\\
(\beta-2\eta)[2+(\beta-2\eta)(b_2+2b_3-b_1)+c]-2\alpha=0,\\
\beta(\beta-2\eta)(b_1+b_2)-c\beta=0,\\
\alpha[(\beta-2\eta)(b_1+b_2)+c]=0.\\
\end{array}\right.
\end{align}
Solving (2.88), (2.88) has no solutions.
Then
\vskip 0.5 true cm
\begin{thm}
$(G_4,g,J)$ is not the first kind algebraic Ricci soliton associated to the connection $\nabla^1$ .
\end{thm}
\vskip 0.5 true cm
\indent By (2.12),(2.13) and (2.86), we have for $(G_4,\nabla^1)$
\begin{align}
\widetilde{{\rm Ric}}^1\left(\begin{array}{c}
e_1\\
e_2\\
e_3
\end{array}\right)=\left(\begin{array}{ccc}
-[1+(\beta-2\eta)(b_3-b_1)]&0&0\\
0&-[1+(\beta-2\eta)(b_2+b_3)]&\frac{b_1+\beta-\alpha-b_3}{2}\\
0&\frac{\alpha+b_3-b_1-\beta}{2}&0
\end{array}\right)\left(\begin{array}{c}
e_1\\
e_2\\
e_3
\end{array}\right).
\end{align}
If $(G_4,g,J)$ is the second kind algebraic Ricci soliton associated to the connection $\nabla^1$, then
$\widetilde{{\rm Ric}}^1=c{\rm Id}+D$, so
\begin{align}
\left\{\begin{array}{l}
De_1=-[1+(\beta-2\eta)(b_3-b_1)+c]e_1,\\
De_2=-[1+(\beta-2\eta)(b_2+b_3)+c]e_2+\frac{b_1+\beta-\alpha-b_3}{2}e_3,\\
De_3=\frac{\alpha+b_3-b_1-\beta}{2}e_2-ce_3.\\
\end{array}\right.
\end{align}
By (2.15) and (2.90), we get
\begin{align}
\left\{\begin{array}{l}
(2\eta-\beta)[\frac{\alpha-\beta}{2}+\frac{3}{2}(b_3-b_1)]-1-c-\frac{\alpha\beta}{2}=0,\\
(\beta-2\eta)[2+(\beta-2\eta)(b_2+2b_3-b_1)+c]-\alpha=0,\\
\beta(\beta-2\eta)(b_1+b_2)+\alpha-c\beta=0,\\
\alpha[(\beta-2\eta)(b_1+b_2)+c]=0.\\
\end{array}\right.
\end{align}
Solving (2.91), (2.91) has no solutions.
Then
\vskip 0.5 true cm
\begin{thm}
$(G_4,g,J)$ is not the second kind algebraic Ricci soliton associated to the connection $\nabla^1$ .
\end{thm}
\vskip 0.5 true cm

\section{ Algebraic Ricci solitons associated to canonical connections and Kobayashi-Nomizu connections on three-dimensional non-unimodular Lorentzian Lie groups}
\vskip 0.5 true cm
\noindent{\bf 3.1 Algebraic Ricci solitons of $G_5$}\\
\vskip 0.5 true cm
By (2.5) and Lemma 4.1 in \cite{BO}, we have for $G_5$, there exists a pseudo-orthonormal basis $\{e_1,e_2,e_3\}$ with $e_3$ timelike such that the Lie
algebra of $G_5$ satisfies
\begin{equation}
[e_1,e_2]=0,~~[e_1,e_3]=\alpha e_1+\beta e_2,~~[e_2,e_3]=\gamma e_1+\delta e_2,~~\alpha+\delta\neq 0,~~\alpha\gamma+\beta\delta=0.
\end{equation}
\vskip 0.5 true cm
\begin{lem}(\cite{Ca},\cite{BO})
The Levi-Civita connection $\nabla$ of $G_5$ is given by
\begin{align}
&\nabla_{e_1}e_1=\alpha e_3,~~\nabla_{e_2}e_1=\frac{\beta+\gamma}{2}e_3,~~\nabla_{e_3}e_1=-\frac{\beta-\gamma}{2}e_2,\\\notag
&\nabla_{e_1}e_2=\frac{\beta+\gamma}{2} e_3,~~\nabla_{e_2}e_2=\delta e_3,~~\nabla_{e_3}e_2=\frac{\beta-\gamma}{2}e_1,\\\notag
&\nabla_{e_1}e_3=\alpha e_1+\frac{\beta+\gamma}{2}  e_2,~~\nabla_{e_2}e_3=\frac{\beta+\gamma}{2} e_1+\delta e_2,~~\nabla_{e_3}e_3=0.
\notag
\end{align}
\end{lem}
By (2.4) and (3.2), we have
\vskip 0.5 true cm
\begin{lem}
For $G_5$, the following equalities hold
\begin{align}
&\nabla_{e_1}(J)e_1=2\alpha e_3,~~\nabla_{e_1}(J)e_2=(\beta+\gamma)e_3,~~\nabla_{e_1}(J)e_3=-2\alpha e_1-(\beta+\gamma)e_2,~~\\\notag
&\nabla_{e_2}(J)e_1=(\beta+\gamma)e_3,~~\nabla_{e_2}(J)e_2=2\delta e_3,~~\nabla_{e_2}(J)e_3=-(\beta+\gamma)e_1-2\delta e_2,~~\\\notag
&\nabla_{e_3}(J)e_1=0,~~\nabla_{e_3}(J)e_2=0,~~\nabla_{e_3}(J)e_3=0.~~\notag
\end{align}
\end{lem}
\vskip 0.5 true cm
By (2.4),(2.5), Lemma 3.1 and Lemma 3.2, we have
\vskip 0.5 true cm
\begin{lem}
The canonical connection $\nabla^0$ of $(G_5,J)$ is given by
\begin{align}
&\nabla^0_{e_1}e_1=0,~~\nabla^0_{e_1}e_2=0,~~\nabla^0_{e_1}e_3=0,\\\notag
&\nabla^0_{e_2}e_1=0,~~\nabla^0_{e_2}e_2=0,~~\nabla^0_{e_2}e_3=0,\\\notag
&\nabla^0_{e_3}e_1=-\frac{\beta-\gamma}{2}e_2,~~\nabla^0_{e_3}e_2=\frac{\beta-\gamma}{2}e_1,~~\nabla^0_{e_3}e_3=0.
\notag
\end{align}
\end{lem}
\vskip 0.5 true cm

\indent By (2.7) and Lemma 3.3, we have
\vskip 0.5 true cm
\begin{lem}
The curvature $R^0$ of the canonical connection $\nabla^0$ of $(G_5,J)$ is flat, that is $R^0(e_i,e_j)e_k=0$.
\end{lem}
\vskip 0.5 true cm
So we get for $(G_5,\nabla^0)$
\begin{align}
{\rm Ric}^0\left(\begin{array}{c}
e_1\\
e_2\\
e_3
\end{array}\right)=\left(\begin{array}{ccc}
0&0&0\\
0&0&0\\
0&0&0
\end{array}\right)\left(\begin{array}{c}
e_1\\
e_2\\
e_3
\end{array}\right).
\end{align}
If $(G_5,g,J)$ is the first kind algebraic Ricci soliton associated to the connection $\nabla^0$, then
${\rm Ric}^0=c{\rm Id}+D$, so
\begin{align}
\left\{\begin{array}{l}
De_1=-ce_1,\\
De_2=-ce_2,\\
De_3=-ce_3.\\
\end{array}\right.
\end{align}
By (2.15) and (3.6), we get
\vskip 0.5 true cm
\begin{thm}
$(G_5,g,J)$ is the first kind algebraic Ricci soliton associated to the connection $\nabla^0$ if and only if $c=0$.
In particular,
\begin{align}
{D}\left(\begin{array}{c}
e_1\\
e_2\\
e_3
\end{array}\right)=\left(\begin{array}{ccc}
0&0&0\\
0&0&0\\
0&0&0
\end{array}\right)\left(\begin{array}{c}
e_1\\
e_2\\
e_3
\end{array}\right).
\end{align}
\end{thm}
\vskip 0.5 true cm
\indent By (2.6), Lemma 3.2 and Lemma 3.3, we have
\vskip 0.5 true cm
\begin{lem}
The Kobayashi-Nomizu connection $\nabla^1$ of $(G_5,J)$ is given by
\begin{align}
&\nabla^1_{e_1}e_1=0,~~\nabla^1_{e_1}e_2=0,~~\nabla^1_{e_1}e_3=0,\\\notag
&\nabla^1_{e_2}e_1=0,~~\nabla^1_{e_2}e_2=0,~~\nabla^1_{e_2}e_3=0,\\\notag
&\nabla^0_{e_3}e_1=-\alpha e_1-\beta e_2,~~\nabla^1_{e_3}e_2=-\gamma e_1-\delta e_2,~~\nabla^1_{e_3}e_3=0.
\notag
\end{align}
\end{lem}
\vskip 0.5 true cm
\begin{lem}
The curvature $R^1$ of the Kobayashi-Nomizu connection $\nabla^1$ of $(G_5,J)$ is flat, that is $R^1(e_i,e_j)e_k=0$.
\end{lem}
\vskip 0.5 true cm
So similar to Theorem 3.5, we have
\vskip 0.5 true cm
\begin{thm}
$(G_5,g,J)$ is the first kind algebraic Ricci soliton associated to the connection $\nabla^1$ if and only if $c=0$.
In particular,
\begin{align}
{D}\left(\begin{array}{c}
e_1\\
e_2\\
e_3
\end{array}\right)=\left(\begin{array}{ccc}
0&0&0\\
0&0&0\\
0&0&0
\end{array}\right)\left(\begin{array}{c}
e_1\\
e_2\\
e_3
\end{array}\right).
\end{align}
\end{thm}

\vskip 0.5 true cm
\noindent{\bf 3.2 Algebraic Ricci solitons of $G_6$}\\
\vskip 0.5 true cm
By (2.6) and Lemma 4.3 in \cite{BO}, we have for $G_6$, there exists a pseudo-orthonormal basis $\{e_1,e_2,e_3\}$ with $e_3$ timelike such that the Lie
algebra of $G_6$ satisfies
\begin{equation}
[e_1,e_2]=\alpha e_2+\beta e_3,~~[e_1,e_3]=\gamma e_2+\delta e_3,~~[e_2,e_3]=0,~~\alpha+\delta\neq 0£¬~~\alpha\gamma-\beta\delta=0.
\end{equation}
\vskip 0.5 true cm
\begin{lem}(\cite{Ca},\cite{BO})
The Levi-Civita connection $\nabla$ of $G_6$ is given by
\begin{align}
&\nabla_{e_1}e_1=0,~~\nabla_{e_2}e_1=-\alpha e_2-\frac{\beta-\gamma}{2}e_3,~~\nabla_{e_3}e_1=\frac{\beta-\gamma}{2}e_2-\delta e_3,\\\notag
&\nabla_{e_1}e_2=\frac{\beta+\gamma}{2} e_3,~~\nabla_{e_2}e_2=\alpha e_1,~~\nabla_{e_3}e_2=-\frac{\beta-\gamma}{2}e_1,\\\notag
&\nabla_{e_1}e_3=\frac{\beta+\gamma}{2}  e_2,~~\nabla_{e_2}e_3=-\frac{\beta-\gamma}{2} e_1,~~\nabla_{e_3}e_3=-\delta e_1.
\notag
\end{align}
\end{lem}
\vskip 0.5 true cm
By (2.4) and (3.11), we have
\vskip 0.5 true cm
\begin{lem}
For $G_6$, the following equalities hold
\begin{align}
&\nabla_{e_1}(J)e_1=0,~~\nabla_{e_1}(J)e_2=(\beta+\gamma)e_3,~~\nabla_{e_1}(J)e_3=-(\beta+\gamma)e_2,~~\\\notag
&\nabla_{e_2}(J)e_1=-(\beta-\gamma) e_3,~~\nabla_{e_2}(J)e_2=0,~~\nabla_{e_2}(J)e_3=(\beta-\gamma) e_1,~~\\\notag
&\nabla_{e_3}(J)e_1=-2\delta e_3,~~\nabla_{e_3}(J)e_2=0,~~\nabla_{e_3}(J)e_3=2\delta e_1.~~\notag
\end{align}
\end{lem}
\vskip 0.5 true cm
By (2.4),(2.5), Lemma 3.9 and Lemma 3.10, we have
\vskip 0.5 true cm
\begin{lem}
The canonical connection $\nabla^0$ of $(G_6,J)$ is given by
\begin{align}
&\nabla^0_{e_1}e_1=0,~~\nabla^0_{e_1}e_2=0,~~\nabla^0_{e_1}e_3=0,\\\notag
&\nabla^0_{e_2}e_1=-\alpha e_2,~~\nabla^0_{e_2}e_2=\alpha e_1,~~\nabla^0_{e_2}e_3=0,\\\notag
&\nabla^0_{e_3}e_1=\frac{\beta-\gamma}{2}e_2,~~\nabla^0_{e_3}e_2=-\frac{\beta-\gamma}{2}e_1,~~\nabla^0_{e_3}e_3=0.
\notag
\end{align}
\end{lem}
\vskip 0.5 true cm
\indent By (2.7) and Lemma 3.11, we have
\vskip 0.5 true cm
\begin{lem}
The curvature $R^0$ of the canonical connection $\nabla^0$ of $(G_6,J)$ is given by
\begin{align}
&R^0(e_1,e_2)e_1=[\alpha^2-\frac{1}{2}\beta(\beta-\gamma)]e_2,~~R^0(e_1,e_2)e_2=[-\alpha^2+\frac{1}{2}\beta(\beta-\gamma)]e_1,~~R^0(e_1,e_2)e_3=0,\\\notag
&R^0(e_1,e_3)e_1=[\gamma\alpha-\frac{1}{2}\delta(\beta-\gamma)] e_2,~~R^0(e_1,e_3)e_2=[-\gamma\alpha+\frac{1}{2}\delta(\beta-\gamma)]  e_1,~~R^0(e_1,e_3)e_3=0,\\\notag
&R^0(e_2,e_3)e_1=0,~~R^0(e_2,e_3)e_2=0,~~R^0(e_2,e_3)e_3=0.\notag
\notag
\end{align}
\end{lem}
\vskip 0.5 true cm
By (2.9),(2.11) and Lemma 3.12, we get for $(G_6,\nabla^0)$
\begin{align}
{\rm Ric}^0\left(\begin{array}{c}
e_1\\
e_2\\
e_3
\end{array}\right)=\left(\begin{array}{ccc}
\frac{1}{2}\beta(\beta-\gamma)-\alpha^2&0&0\\
0&\frac{1}{2}\beta(\beta-\gamma)-\alpha^2&0\\
0&-\gamma\alpha+\frac{1}{2}\delta(\beta-\gamma)&0
\end{array}\right)\left(\begin{array}{c}
e_1\\
e_2\\
e_3
\end{array}\right).
\end{align}
If $(G_6,g,J)$ is the first kind algebraic Ricci soliton associated to the connection $\nabla^0$, then
${\rm Ric}^0=c{\rm Id}+D$, so
\begin{align}
\left\{\begin{array}{l}
De_1=[\frac{1}{2}\beta(\beta-\gamma)-\alpha^2-c]e_1,\\
De_2=[\frac{1}{2}\beta(\beta-\gamma)-\alpha^2-c]e_2,\\
De_3=[-\gamma\alpha+\frac{1}{2}\delta(\beta-\gamma)]e_2-ce_3.\\
\end{array}\right.
\end{align}
By (2.15) and (3.16), we get
\begin{align}
\left\{\begin{array}{l}
\alpha[\frac{1}{2}\beta(\beta-\gamma)-\alpha^2-c]-\beta[-\gamma\alpha+\frac{1}{2}\delta(\beta-\gamma)]=0,\\
\beta[\beta(\beta-\gamma)-2\alpha^2-c]=0,\\
-c\gamma+[-\gamma\alpha+\frac{1}{2}\delta(\beta-\gamma)](\alpha-\delta)=0,\\
\delta[\frac{1}{2}\beta(\beta-\gamma)-\alpha^2-c]+\beta[-\gamma\alpha+\frac{1}{2}\delta(\beta-\gamma)]=0.\\
\end{array}\right.
\end{align}
Solving (3.17), then
\vskip 0.5 true cm
\begin{thm}
$(G_6,g,J)$ is the first kind algebraic Ricci soliton associated to the connection $\nabla^0$ if and only if \\
(i)$\alpha=\beta=\gamma=c=0$, $\delta\neq 0$,
in particular,
 \begin{align}
{D}\left(\begin{array}{c}
e_1\\
e_2\\
e_3
\end{array}\right)=\left(\begin{array}{ccc}
0&0&0\\
0&0&0\\
0&0&0
\end{array}\right)\left(\begin{array}{c}
e_1\\
e_2\\
e_3 \notag \\
\end{array}\right).
\end{align}
(ii)$\alpha\neq 0$, $\beta=\gamma=0$, $\alpha^2+c=0$, $\alpha+\delta\neq 0$,
in particular,
 \begin{align}
{D}\left(\begin{array}{c}
e_1\\
e_2\\
e_3
\end{array}\right)=\left(\begin{array}{ccc}
0&0&0\\
0&0&0\\
0&0&\alpha^2
\end{array}\right)\left(\begin{array}{c}
e_1\\
e_2\\
e_3 \notag \\
\end{array}\right).
\end{align}
(iii)$\alpha\neq 0$, $\beta\neq 0$, $\gamma=\delta=0$, $\beta^2=2\alpha^2$, $c=0$,
in particular,
 \begin{align}
{D}\left(\begin{array}{c}
e_1\\
e_2\\
e_3
\end{array}\right)=\left(\begin{array}{ccc}
0&0&0\\
0&0&0\\
0&0&0
\end{array}\right)\left(\begin{array}{c}
e_1\\
e_2\\
e_3 \notag \\
\end{array}\right).
\end{align}
\end{thm}
\vskip 0.5 true cm
\indent By (2.12),(2.13) and (3.15), we have for $(G_6,\nabla^0)$
\begin{align}
\widetilde{{\rm Ric}}^0\left(\begin{array}{c}
e_1\\
e_2\\
e_3
\end{array}\right)==\left(\begin{array}{ccc}
\frac{1}{2}\beta(\beta-\gamma)-\alpha^2&0&0\\
0&\frac{1}{2}\beta(\beta-\gamma)-\alpha^2&\frac{1}{2}[\gamma\alpha-\frac{1}{2}\delta(\beta-\gamma)]\\
0&\frac{1}{2}[-\gamma\alpha+\frac{1}{2}\delta(\beta-\gamma)]&0
\end{array}\right)\left(\begin{array}{c}
e_1\\
e_2\\
e_3
\end{array}\right).
\end{align}
If $(G_6,g,J)$ is the second kind algebraic Ricci soliton associated to the connection $\nabla^0$, then
$\widetilde{{\rm Ric}}^0=c{\rm Id}+D$, so
\begin{align}
\left\{\begin{array}{l}
De_1=[\frac{1}{2}\beta(\beta-\gamma)-\alpha^2-c]e_1,\\
De_2=[\frac{1}{2}\beta(\beta-\gamma)-\alpha^2-c]e_2+\frac{1}{2}[\gamma\alpha-\frac{1}{2}\delta(\beta-\gamma)]e_3,\\
De_3=\frac{1}{2}[-\gamma\alpha+\frac{1}{2}\delta(\beta-\gamma)]e_2-ce_3.\\
\end{array}\right.
\end{align}
By (2.15) and (3.19), we get
\begin{align}
\left\{\begin{array}{l}
\alpha[\frac{1}{2}\beta(\beta-\gamma)-\alpha^2-c]+\frac{1}{2}(\beta+\gamma)[\gamma\alpha-\frac{1}{2}\delta(\beta-\gamma)]=0,\\
\beta[\beta(\beta-\gamma)-2\alpha^2-c]+\frac{1}{2}(\delta-\alpha)[\gamma\alpha-\frac{1}{2}\delta(\beta-\gamma)]=0,\\
-c\gamma+\frac{1}{2}(\alpha-\delta)[-\gamma\alpha+\frac{1}{2}\delta(\beta-\gamma)]=0,\\
\delta[\frac{1}{2}\beta(\beta-\gamma)-\alpha^2-c]+\frac{1}{2}(\beta+\gamma)[-\gamma\alpha+\frac{1}{2}\delta(\beta-\gamma)]=0.\\
\end{array}\right.
\end{align}
Solving (3.20), we get
\vskip 0.5 true cm
\begin{thm}
$(G_6,g,J)$ is the second kind algebraic Ricci soliton associated to the connection $\nabla^0$ if and only if \\
(i)$\beta=\gamma=0$, $\alpha^2+c=0$, $\alpha+\delta\neq 0$,
in particular,
 \begin{align}
{D}\left(\begin{array}{c}
e_1\\
e_2\\
e_3
\end{array}\right)=\left(\begin{array}{ccc}
0&0&0\\
0&0&0\\
0&0&-c
\end{array}\right)\left(\begin{array}{c}
e_1\\
e_2\\
e_3 \notag \\
\end{array}\right).
\end{align}
(ii)$\gamma=\delta=c=0$, $\alpha\neq 0$, $\beta\neq 0$, $\beta^2=2\alpha^2$,
in particular,
 \begin{align}
{D}\left(\begin{array}{c}
e_1\\
e_2\\
e_3
\end{array}\right)=\left(\begin{array}{ccc}
0&0&0\\
0&0&0\\
0&0&0
\end{array}\right)\left(\begin{array}{c}
e_1\\
e_2\\
e_3 \notag \\
\end{array}\right).
\end{align}
\end{thm}
\vskip 0.5 true cm
\indent By (2.6), Lemma 3.10 and Lemma 3.11, we have
\vskip 0.5 true cm
\begin{lem}
The Kobayashi-Nomizu connection $\nabla^1$ of $(G_6,J)$ is given by
\begin{align}
&\nabla^1_{e_1}e_1=0,~~\nabla^1_{e_1}e_2=0,~~\nabla^1_{e_1}e_3=\delta e_3,\\\notag
&\nabla^1_{e_2}e_1=-\alpha e_2,~~\nabla^1_{e_2}e_2=\alpha e_1,~~\nabla^1_{e_2}e_3=0,\\\notag
&\nabla^1_{e_3}e_1=-\gamma e_2,~~\nabla^1_{e_3}e_2=0,~~\nabla^1_{e_3}e_3=0.
\notag
\end{align}
\end{lem}
\vskip 0.5 true cm
\indent By (2.8) and Lemma 3.15, we have
\vskip 0.5 true cm
\begin{lem}
The curvature $R^1$ of the Kobayashi-Nomizu connection $\nabla^1$ of $(G_6,J)$ is given by
\begin{align}
&R^1(e_1,e_2)e_1=(\alpha^2+\beta\gamma)e_2,~~R^1(e_1,e_2)e_2=-\alpha^2e_1,~~R^1(e_1,e_2)e_3=0,\\\notag
&R^1(e_1,e_3)e_1=\gamma(\alpha+\delta)e_2,~~R^1(e_1,e_3)e_2=-\alpha\gamma e_1,~~R^1(e_1,e_3)e_3=0,\\\notag
&R^1(e_2,e_3)e_1=-\alpha\gamma e_1,~~R^1(e_2,e_3)e_2=\alpha\gamma e_2,~~R^1(e_2,e_3)e_3=0.\notag
\notag
\end{align}
\end{lem}
\vskip 0.5 true cm
By (2.10),(2.11) and Lemma 3.16, we get for $(G_6,\nabla^1)$
\begin{align}
{\rm Ric}^1\left(\begin{array}{c}
e_1\\
e_2\\
e_3
\end{array}\right)=\left(\begin{array}{ccc}
-(\alpha^2+\beta\gamma)&0&0\\
0&-\alpha^2&0\\
0&0&0
\end{array}\right)\left(\begin{array}{c}
e_1\\
e_2\\
e_3
\end{array}\right).
\end{align}
If $(G_6,g,J)$ is the first kind algebraic Ricci soliton associated to the connection $\nabla^1$, then
${\rm Ric}^1=c{\rm Id}+D$, so
\begin{align}
\left\{\begin{array}{l}
De_1=-(\alpha^2+\beta\gamma+c)e_1,\\
De_2=-(\alpha^2+c) e_2,\\
De_3=-ce_3.\\
\end{array}\right.
\end{align}

By (2.15) and (3.24), we get
\begin{align}
\left\{\begin{array}{l}
\alpha(\alpha^2+\beta\gamma+c)=0,\\
\beta(2\alpha^2+\beta\gamma+c)=0,\\
(\beta\gamma+c)\gamma=0,\\
(\alpha^2+\beta\gamma+c)\delta=0.\\
\end{array}\right.
\end{align}
\noindent Solving (3.25), we get
\vskip 0.5 true cm
\begin{thm}
$(G_6,g,J)$ is the first kind algebraic Ricci soliton associated to the connection $\nabla^1$ if and only if \\
(i)$\alpha=\beta=c=0$, $\delta\neq 0$,
in particular,
 \begin{align}
{D}\left(\begin{array}{c}
e_1\\
e_2\\
e_3
\end{array}\right)=\left(\begin{array}{ccc}
0&0&0\\
0&0&0\\
0&0&0
\end{array}\right)\left(\begin{array}{c}
e_1\\
e_2\\
e_3 \notag \\
\end{array}\right).
\end{align}
(ii)$\alpha\neq 0$, $\beta=\gamma=0$, $\alpha^2+c=0$, $\alpha+\delta\neq 0$,
in particular,
 \begin{align}
{D}\left(\begin{array}{c}
e_1\\
e_2\\
e_3
\end{array}\right)=\left(\begin{array}{ccc}
0&0&0\\
0&0&0\\
0&0&\alpha^2
\end{array}\right)\left(\begin{array}{c}
e_1\\
e_2\\
e_3 \notag \\
\end{array}\right).
\end{align}
\end{thm}
\vskip 0.5 true cm

\noindent{\bf 3.3 Algebraic Ricci solitons of $G_7$}\\
\vskip 0.5 true cm
By (2.7) and Lemma 4.5 in \cite{BO}, we have for $G_7$, there exists a pseudo-orthonormal basis $\{e_1,e_2,e_3\}$ with $e_3$ timelike such that the Lie
algebra of $G_7$ satisfies
\begin{equation}
[e_1,e_2]=-\alpha e_1-\beta e_2-\beta e_3,~~[e_1,e_3]=\alpha e_1+\beta e_2+\beta e_3,~~[e_2,e_3]=\gamma e_1+\delta e_2+\delta e_3,,~~\alpha+\delta\neq 0,~~\alpha\gamma=0.
\end{equation}

\vskip 0.5 true cm
\begin{lem}(\cite{Ca},\cite{BO})
The Levi-Civita connection $\nabla$ of $G_7$ is given by
\begin{align}
&\nabla_{e_1}e_1=\alpha e_2+\alpha e_3,~~\nabla_{e_2}e_1=\beta e_2+(\beta+\frac{\gamma}{2})e_3,~~\nabla_{e_3}e_1=-(\beta-\frac{\gamma}{2})e_2-\beta e_3,\\\notag
&\nabla_{e_1}e_2=-\alpha e_1+\frac{\gamma}{2} e_3,~~\nabla_{e_2}e_2=-\beta e_1+\delta e_3,~~\nabla_{e_3}e_2=(\beta-\frac{\gamma}{2})e_1-\delta e_3,\\\notag
&\nabla_{e_1}e_3=\alpha e_1+\frac{\gamma}{2} e_2,~~\nabla_{e_2}e_3=(\beta+\frac{\gamma}{2})e_1+\delta e_2,~~\nabla_{e_3}e_3=-\beta e_1-\delta e_2.
\notag
\end{align}
\end{lem}
\vskip 0.5 true cm
By (2.4) and (3.27), we have
\vskip 0.5 true cm
\begin{lem}
For $G_7$, the following equalities hold
\begin{align}
&\nabla_{e_1}(J)e_1=2\alpha e_3,~~\nabla_{e_1}(J)e_2=\gamma e_3,~~\nabla_{e_1}(J)e_3=-2\alpha e_1-\gamma e_2,~~\\\notag
&\nabla_{e_2}(J)e_1=2(\beta+\frac{\gamma}{2}) e_3,~~\nabla_{e_2}(J)e_2=2\delta e_3,~~\nabla_{e_2}(J)e_3=-2(\beta+\frac{\gamma}{2})e_1-2\delta e_2,~~\\\notag
&\nabla_{e_3}(J)e_1=-2\beta e_3,~~\nabla_{e_3}(J)e_2=-2\delta e_3, ,~~\nabla_{e_3}(J)e_3=2\beta e_1+2\delta e_2.~~\notag
\end{align}
\end{lem}
\vskip 0.5 true cm
By (2.4),(2.5), Lemma 3.18 and Lemma 3.19, we have
\vskip 0.5 true cm
\begin{lem}
The canonical connection $\nabla^0$ of $(G_7,J)$ is given by
\begin{align}
&\nabla^0_{e_1}e_1=\alpha e_2,~~\nabla^0_{e_1}e_2=-\alpha e_1,~~\nabla^0_{e_1}e_3=0,\\\notag
&\nabla^0_{e_2}e_1=\beta e_2,~~\nabla^0_{e_2}e_2=-\beta e_1,~~\nabla^0_{e_2}e_3=0,\\\notag
&\nabla^0_{e_3}e_1=-(\beta-\frac{\gamma}{2})e_2,~~\nabla^0_{e_3}e_2=(\beta-\frac{\gamma}{2})e_1,~~\nabla^0_{e_3}e_3=0.
\notag
\end{align}
\end{lem}
\vskip 0.5 true cm
\indent By (2.7) and Lemma 3.20, we have
\vskip 0.5 true cm
\begin{lem}
The curvature $R^0$ of the canonical connection $\nabla^0$ of $(G_7,J)$ is given by
\begin{align}
&R^0(e_1,e_2)e_1=(\alpha^2+\frac{\beta\gamma}{2})e_2,~~R^0(e_1,e_2)e_2=-(\alpha^2+\frac{\beta\gamma}{2})e_1,~~R^0(e_1,e_2)e_3=0,\\\notag
&R^0(e_1,e_3)e_1=-(\alpha^2+\frac{\beta\gamma}{2}) e_2,~~R^0(e_1,e_3)e_2=(\alpha^2+\frac{\beta\gamma}{2}) e_1,~~R^0(e_1,e_3)e_3=0,\\\notag
&R^0(e_2,e_3)e_1=-(\gamma\alpha+\frac{\delta\gamma}{2})e_2,~~R^0(e_2,e_3)e_2=(\gamma\alpha+\frac{\delta\gamma}{2})e_1,~~R^0(e_2,e_3)e_3=0.\notag
\notag
\end{align}
\end{lem}
\vskip 0.5 true cm
By (2.9),(2.11) and Lemma 3.21, we get for $(G_7,\nabla^0)$
\begin{align}
{\rm Ric}^0\left(\begin{array}{c}
e_1\\
e_2\\
e_3
\end{array}\right)=\left(\begin{array}{ccc}
-(\alpha^2+\frac{\beta\gamma}{2})&0&0\\
0&-(\alpha^2+\frac{\beta\gamma}{2})&0\\
-(\gamma\alpha+\frac{\delta\gamma}{2})&(\alpha^2+\frac{\beta\gamma}{2})&0
\end{array}\right)\left(\begin{array}{c}
e_1\\
e_2\\
e_3
\end{array}\right).
\end{align}
\noindent If $(G_7,g,J)$ is the first kind algebraic Ricci soliton associated to the connection $\nabla^0$, then
${\rm Ric}^0=c{\rm Id}+D$, so
\begin{align}
\left\{\begin{array}{l}
De_1=-(\alpha^2+\frac{\beta\gamma}{2}+c)e_1,\\
De_2=-(\alpha^2+\frac{\beta\gamma}{2}+c)e_2,\\
De_3=-(\gamma\alpha+\frac{\delta\gamma}{2})e_1+(\alpha^2+\frac{\beta\gamma}{2})e_2-ce_3.\\
\end{array}\right.
\end{align}
By (2.15) and (3.32), we get
\begin{align}
\left\{\begin{array}{l}
\alpha(\alpha^2+\frac{\beta\gamma}{2}+c)-\beta(\gamma\alpha+\frac{\delta\gamma}{2})=0,\\
\beta(2\alpha^2+\beta\gamma+c)=0,\\
c\gamma+(\alpha-\delta)(\gamma\alpha+\frac{\delta\gamma}{2})=0.\\
(\alpha^2+\frac{\beta\gamma}{2}+c)\delta+\beta(\gamma\alpha+\frac{\delta\gamma}{2})=0
\end{array}\right.
\end{align}
Solving (3.33), we get
\vskip 0.5 true cm
\begin{thm}
$(G_7,g,J)$ is the first kind algebraic Ricci soliton associated to the connection $\nabla^0$ if and only if\\
 (i)$\alpha=\gamma=c=0$, $\beta\neq 0$, $\delta\neq 0$,
in particular,
 \begin{align}
{D}\left(\begin{array}{c}
e_1\\
e_2\\
e_3
\end{array}\right)=\left(\begin{array}{ccc}
0&0&0\\
0&0&0\\
0&0&0
\end{array}\right)\left(\begin{array}{c}
e_1\\
e_2\\
e_3 \notag \\
\end{array}\right).
\end{align}
(ii)$\beta=\gamma=0$, $\alpha^2+c=0$, $\alpha+\delta\neq 0$,
in particular,
 \begin{align}
{D}\left(\begin{array}{c}
e_1\\
e_2\\
e_3
\end{array}\right)=\left(\begin{array}{ccc}
0&0&0\\
0&0&0\\
0&\alpha^2&\alpha^2
\end{array}\right)\left(\begin{array}{c}
e_1\\
e_2\\
e_3 \notag \\
\end{array}\right).
\end{align}
\end{thm}
\vskip 0.5 true cm
\indent By (2.12),(2.13) and (3.31), we have for $(G_7,\nabla^0)$
\begin{align}
\widetilde{{\rm Ric}}^0\left(\begin{array}{c}
e_1\\
e_2\\
e_3
\end{array}\right)=\left(\begin{array}{ccc}
-(\alpha^2+\frac{\beta\gamma}{2})&0&\frac{1}{2}(\gamma\alpha+\frac{\delta\gamma}{2})\\
0&-(\alpha^2+\frac{\beta\gamma}{2})&-\frac{1}{2}(\alpha^2+\frac{\beta\gamma}{2})\\
-\frac{1}{2}(\gamma\alpha+\frac{\delta\gamma}{2})&\frac{1}{2}(\alpha^2+\frac{\beta\gamma}{2})&0
\end{array}\right)\left(\begin{array}{c}
e_1\\
e_2\\
e_3
\end{array}\right).
\end{align}
\noindent If $(G_7,g,J)$ is the second kind algebraic Ricci soliton associated to the connection $\nabla^0$, then
$\widetilde{{\rm Ric}}^0=c{\rm Id}+D$, so
\begin{align}
\left\{\begin{array}{l}
De_1=-(\alpha^2+\frac{\beta\gamma}{2}+c)e_1+\frac{1}{2}(\gamma\alpha+\frac{\delta\gamma}{2})e_3,\\
De_2=-(\alpha^2+\frac{\beta\gamma}{2}+c)e_2-\frac{1}{2}(\alpha^2+\frac{\beta\gamma}{2})e_3,\\
De_3=-\frac{1}{2}(\gamma\alpha+\frac{\delta\gamma}{2})e_1+\frac{1}{2}(\alpha^2+\frac{\beta\gamma}{2})e_2-ce_3.\\
\end{array}\right.
\end{align}
By (2.15) and (3.35), we get
\begin{align}
\left\{\begin{array}{l}
\alpha(\frac{\alpha^2}{2}+\frac{\beta\gamma}{4}+c)-\frac{1}{2}(\beta+\gamma)(\gamma\alpha+\frac{\delta\gamma}{2})=0,\\
\beta(\alpha^2+\frac{\beta\gamma}{2}+c)-\frac{1}{2}\delta(\gamma\alpha+\frac{\delta\gamma}{2})=0,\\
\alpha(\frac{\alpha^2}{2}+\frac{\beta\gamma}{4}+c)-\frac{1}{2}\beta(\gamma\alpha+\frac{\delta\gamma}{2})=0,\\
\beta(\alpha^2+\frac{\beta\gamma}{2}+c)=0,\\
\gamma c-\frac{1}{4}\gamma\delta^2=0,\\
\delta(\frac{\alpha^2}{2}+\frac{\beta\gamma}{4}+c)+\frac{1}{2}\beta(\gamma\alpha+\frac{\delta\gamma}{2})=0,\\
\delta(\frac{\alpha^2}{2}+\frac{\beta\gamma}{4}+c)+\frac{1}{2}(\beta+\gamma)(\gamma\alpha+\frac{\delta\gamma}{2})=0,\\
\end{array}\right.
\end{align}
Solving (3.36), we get
\vskip 0.5 true cm
\begin{thm}
$(G_7,g,J)$ is the second kind algebraic Ricci soliton associated to the connection $\nabla^0$ if and only if \\
(i)$\alpha=\gamma=c=0$, $\delta\neq 0$,
in particular,
 \begin{align}
{D}\left(\begin{array}{c}
e_1\\
e_2\\
e_3
\end{array}\right)=\left(\begin{array}{ccc}
0&0&0\\
0&0&0\\
0&0&0
\end{array}\right)\left(\begin{array}{c}
e_1\\
e_2\\
e_3 \notag \\
\end{array}\right).
\end{align}
(ii)$\alpha\neq 0$, $\beta=\gamma=0$, $\frac{\alpha^2}{2}+c=0$, $\alpha+\delta\neq 0$,
in particular,
 \begin{align}
{D}\left(\begin{array}{c}
e_1\\
e_2\\
e_3
\end{array}\right)=\left(\begin{array}{ccc}
-\frac{\alpha^2}{2}&0&0\\
0&-\frac{\alpha^2}{2}&-\frac{\alpha^2}{2}\\
0&\frac{\alpha^2}{2}&\frac{\alpha^2}{2}
\end{array}\right)\left(\begin{array}{c}
e_1\\
e_2\\
e_3 \notag \\
\end{array}\right).
\end{align}
\end{thm}
\vskip 0.5 true cm
\indent By (2.6), Lemma 3.19 and Lemma 3.20, we have
\vskip 0.5 true cm
\begin{lem}
The Kobayashi-Nomizu connection $\nabla^1$ of $(G_7,J)$ is given by
\begin{align}
&\nabla^1_{e_1}e_1=\alpha e_2,~~\nabla^1_{e_1}e_2=-\alpha e_1,~~\nabla^1_{e_1}e_3=\beta e_3,\\\notag
&\nabla^1_{e_2}e_1=\beta e_2,~~\nabla^1_{e_2}e_2=-\beta e_1,~~\nabla^1_{e_2}e_3=\delta e_3,\\\notag
&\nabla^1_{e_3}e_1=-\alpha e_1-\beta e_2,~~\nabla^1_{e_3}e_2=-\gamma e_1-\delta e_2,~~\nabla^1_{e_3}e_3=0.
\notag
\end{align}
\end{lem}
\indent By (2.8) and Lemma 3.24, we have
\vskip 0.5 true cm
\begin{lem}
The curvature $R^1$ of the Kobayashi-Nomizu connection $\nabla^1$ of $(G_7,J)$ is given by
\begin{align}
&R^1(e_1,e_2)e_1=-\alpha\beta e_1+\alpha^2 e_2,~~R^1(e_1,e_2)e_2=-(\alpha^2+\beta^2+\beta\gamma) e_1-\beta\delta e_2,\\\notag
&R^1(e_1,e_2)e_3=\beta(\alpha+\delta)e_3,
R^1(e_1,e_3)e_1=(2\alpha\beta+\alpha\gamma)e_1+(-2\alpha^2+\alpha\delta)e_2,\\\notag
&R^1(e_1,e_3)e_2=(\alpha\delta+\beta^2+\beta\gamma)e_1
-(\gamma\alpha+\alpha\beta-\beta\delta)e_2,~~R^1(e_1,e_3)e_3=-\beta(\alpha+\delta)e_3,\\\notag
&R^1(e_2,e_3)e_1=(\beta^2+\beta\gamma+\alpha\delta)e_1-(\alpha\beta-\beta\delta+\alpha\gamma)e_2,\\\notag
&R^1(e_2,e_3)e_2= (2\beta\delta-\alpha\beta+\alpha\gamma+\delta\gamma)e_1-(\beta\gamma+\beta^2-\delta^2)e_2,\\\notag
&R^1(e_2,e_3)e_3=-(\beta\gamma+\delta^2)e_3.\notag
\notag
\end{align}
\end{lem}
\vskip 0.5 true cm
By (2.10),(2.11) and Lemma 3.25, we get for $(G_7,\nabla^1)$
\begin{align}
{\rm Ric}^1\left(\begin{array}{c}
e_1\\
e_2\\
e_3
\end{array}\right)=\left(\begin{array}{ccc}
-\alpha^2&\beta\delta&-\beta(\alpha+\delta)\\
-\alpha\beta&-(\alpha^2+\beta^2+\beta\gamma)&-(\beta\gamma+\delta^2)\\
\beta(\alpha+\delta)&\delta(\alpha+\delta)&0
\end{array}\right)\left(\begin{array}{c}
e_1\\
e_2\\
e_3
\end{array}\right).
\end{align}
\noindent If $(G_7,g,J)$ is the first kind algebraic Ricci soliton associated to the connection $\nabla^1$, then
${\rm Ric}^1=c{\rm Id}+D$, so
\begin{align}
\left\{\begin{array}{l}
De_1=-(\alpha^2+c)e_1+\beta\delta e_2-\beta(\alpha+\delta)e_3,\\
De_2=-\alpha\beta e_1-(\alpha^2+\beta^2+\beta\gamma+c) e_2-(\beta\gamma+\delta^2)e_3,\\
De_3=\beta(\alpha+\delta) e_1+\delta(\alpha+\delta)e_2-ce_3.\\
\end{array}\right.
\end{align}
By (2.15) and (3.40), we get
\begin{align}
\left\{\begin{array}{l}
\alpha(\alpha^2+\beta^2+c)+\beta^2\delta+\beta\gamma\delta-\alpha\delta^2=0,\\
\beta\alpha^2+\beta c+3\alpha\beta\delta+\beta\delta^2-\beta^2\gamma=0,\\
\beta(\alpha^2+\beta^2+c-\delta^2-\beta\gamma)=0,\\
\alpha c+\beta^2\delta-\beta\gamma\delta+\alpha^2 \delta+\alpha\delta^2=0,\\
\beta(\beta^2+\beta\gamma-3\alpha\delta-c-\delta^2)=0,\\
\beta(\beta\gamma+\delta^2-c)=0,\\
\gamma(\alpha^2+\beta^2+\beta\gamma+c)+\beta\delta^2-\alpha\beta\delta=0,\\
\delta(\beta\gamma+c+\delta\alpha+\delta^2-\beta^2)=0,\\
\delta(\alpha^2+\beta^2-\delta^2+c)-(\alpha+\delta)(\beta\gamma+\beta^2)+\alpha\beta^2=0.\\
\end{array}\right.
\end{align}
\noindent Solving (3.41), we get
\vskip 0.5 true cm
\begin{thm}
$(G_7,g,J)$ is the first kind algebraic Ricci soliton associated to the connection $\nabla^1$ if and only if
$\alpha\neq 0$, $\beta=\gamma=0$, $\alpha=2\delta$, $c=-3\delta^2$, in particular
 \begin{align}
{D}\left(\begin{array}{c}
e_1\\
e_2\\
e_3
\end{array}\right)=\left(\begin{array}{ccc}
-\delta^2&0&0\\
0&-\delta^2&-\delta^2\\
0&3\delta^2&3\delta^2
\end{array}\right)\left(\begin{array}{c}
e_1\\
e_2\\
e_3 \notag \\
\end{array}\right).
\end{align}
\end{thm}
\vskip 0.5 true cm
\indent By (2.12),(2.13) and (3.39), we have for $(G_7,\nabla^1)$
\begin{align}
\widetilde{{\rm Ric}}^1\left(\begin{array}{c}
e_1\\
e_2\\
e_3
\end{array}\right)=\left(\begin{array}{ccc}
-\alpha^2&\frac{1}{2}(\beta\delta-\alpha\beta)&-\beta(\alpha+\delta)\\
\frac{1}{2}(\beta\delta-\alpha\beta)&-(\alpha^2+\beta^2+\beta\gamma)&-\frac{1}{2}(\beta\gamma+\alpha\delta+2\delta^2)\\
\beta(\alpha+\delta)&\frac{1}{2}(\beta\gamma+\alpha\delta+2\delta^2)&0
\end{array}\right)\left(\begin{array}{c}
e_1\\
e_2\\
e_3
\end{array}\right).
\end{align}
\noindent If $(G_7,g,J)$ is the second kind algebraic Ricci soliton associated to the connection $\nabla^1$, then
$\widetilde{{\rm Ric}}^1=c{\rm Id}+D$, so
\begin{align}
\left\{\begin{array}{l}
De_1=-(\alpha^2+c)e_1+\frac{1}{2}(\beta\delta-\alpha\beta)e_2-\beta(\alpha+\delta)e_3,\\
De_2=\frac{1}{2}(\beta\delta-\alpha\beta)e_1-(\alpha^2+\beta^2+\beta\gamma+c)e_2-\frac{1}{2}(\beta\gamma+\alpha\delta+2\delta^2)e_3,\\
De_3=\beta(\alpha+\delta)e_1+\frac{1}{2}(\beta\gamma+\alpha\delta+2\delta^2)e_2-ce_3.\\
\end{array}\right.
\end{align}

By (2.15) and (3.43), we get
\begin{align}
\left\{\begin{array}{l}
\alpha(\alpha^2+\beta^2+\beta\gamma+c)+(\beta^2+\beta\gamma)(\alpha+\delta)+\frac{\beta}{2}(\beta\delta-\alpha\beta)-\frac{\alpha}{2}(\beta\gamma+\alpha\delta
+2\delta^2)=0,\\
\beta(\frac{\alpha^2}{2}+c+\frac{3}{2}\alpha\delta+\delta^2)=0,\\
\beta(2\alpha^2+\beta^2+\beta\gamma+c)+\beta(\alpha+\delta)(\delta-\alpha)-\beta(\beta\gamma+\alpha\delta+2\delta^2)=0,\\
\alpha c+\frac{1}{2}(\beta-\gamma)(\beta\delta-\alpha\beta)+\beta^2(\alpha+\delta)+\frac{\alpha}{2}(\beta\gamma+\alpha\delta+2\delta^2)=0,\\
\frac{1}{2}(\alpha-\delta)(\beta\delta-\alpha\beta)-\beta(\beta^2+\beta\gamma-c)+\beta(\beta\gamma+\alpha\delta+2\delta^2)=0,\\
\beta(c-\frac{\alpha\delta}{2}-\frac{\delta^2}{2})=0,\\
\frac{1}{2}(\delta-\alpha)(\beta\delta-\alpha\beta)+\beta(\alpha+\delta)(\delta-\alpha)+\gamma(\beta^2+\beta\gamma+c)=0,\\
\frac{1}{2}(\gamma-\beta)(\beta\delta-\alpha\beta)+\delta c+\frac{\delta}{2}(\beta\gamma+\alpha\delta+2\delta^2)-\beta^2(\alpha+\delta)=0,\\
\frac{1}{2}\beta(\beta\delta-\alpha\beta)+\frac{1}{2}\delta(\beta\gamma+\alpha\delta+2\delta^2)-\delta(\alpha^2+\beta^2+\beta\gamma+c)+(\beta^2+\beta\gamma)
(\alpha+\delta)=0.\\
\end{array}\right.
\end{align}
\noindent Solving (3.44), we get
\vskip 0.5 true cm
\begin{thm}
$(G_7,g,J)$ is the second kind algebraic Ricci soliton associated to the connection $\nabla^1$ if and only if
$\alpha\neq 0$, $\beta=\gamma=0$, $\alpha=2\delta$, $\alpha^2+2c=0$, in particular
\notag\begin{align}
{D}\left(\begin{array}{c}
e_1\\
e_2\\
e_3
\end{array}\right)=\left(\begin{array}{ccc}
-\frac{\alpha^2}{2}&0&0\\
0&-\frac{\alpha^2}{2}&-\frac{\alpha^2}{2}\\
0&\frac{\alpha^2}{2}&\frac{\alpha^2}{2}
\end{array}\right)\left(\begin{array}{c}
e_1\\
e_2\\
e_3
\end{array}\right).
\end{align}
\end{thm}
\vskip 1 true cm

\section{Acknowledgements}

The author was supported in part by NSFC No.11771070.

\vskip 1 true cm


\bigskip
\bigskip

\noindent {\footnotesize {\it Y. Wang} \\
{School of Mathematics and Statistics, Northeast Normal University, Changchun 130024, China}\\
{Email: wangy581@nenu.edu.cn}


\begin{thebibliography}{20}

\bibitem{BO}
W. Batat, K. Onda, {\it Algebraic Ricci solitons of three-dimensional Lorentzian Lie groups}, J. Geom. Phys. 114 (2017), 138-152.

 \bibitem{Ca1} G. Calvaruso, {\it Homogeneous structures on three-dimensional Lorentzian manifolds}, J. Geom. Phys 57 (4) (2007) 1279-1291.

\bibitem{Ca} G. Calvaruso, {\it Einstein-like metrics on three-dimensional homogeneous Lorentzian manifolds}, Geom. Dedicata 127 (2007) 99-119.

\bibitem{CP} L.A. Cordero, P.E. Parker, {\it Left-invariant Lorentzian metrics on 3-dimensional Lie groups}, Rend. Mat. Appl. (7) 17 (1997) 129-155.

\bibitem{ES}
F. Etayo, R. Santamaria, {\it Distinguished connections on $(J^2=\pm 1)$-metric manifolds}, Arch. Math. (Brno) 52 (2016), no. 3, 159-203.

\bibitem{La}
J. Lauret, {\it Ricci soliton homogeneous nilmanifolds}, Math. Ann. 319 (4)(2001), 715-733.


\bibitem{On}
K. Onda, {\it Examples of algebraic Ricci solitons in the pseudo-Riemannian case}, Acta Math. Hungar. 144 (2014), no. 1, 247-265









\end{thebibliography}
\end{document}